\documentclass[12pt]{article}
\usepackage[round,sort,authoryear]{natbib}

\usepackage{amsmath,amssymb,amsthm,amsfonts, fullpage, setspace}

\usepackage{eucal}
\usepackage{subeqnarray,bm}
 \usepackage[utf8]{inputenc}
\usepackage{graphicx}[final]
\usepackage{psfrag} 
\usepackage{color}
\usepackage{pdfpages}

\usepackage{epstopdf}

\usepackage[colorlinks]{hyperref}

\usepackage{hyperref}
\hypersetup{colorlinks=red,citecolor=blue,pdfstartview=FitH, pdfpagemode=None}
\usepackage{multirow}
\usepackage[bottom]{footmisc}
\usepackage{ifthen} 
\usepackage{sectsty}
\usepackage{subfig}
\usepackage{float}
\usepackage[english]{babel}
\usepackage[nottoc]{tocbibind}
\usepackage{calrsfs}
\usepackage{soul}
 \DeclareGraphicsExtensions{.eps,.png,.pdf}
\usepackage{mwe}
\makeatletter
\def\maxwidth{ %
  \ifdim\Gin@nat@width>\linewidth
    \linewidth
  \else
    \Gin@nat@width
  \fi
}
\makeatother

\definecolor{fgcolor}{rgb}{0.345, 0.345, 0.345}

\allsectionsfont{\normalsize\bfseries}

\usepackage{framed}
\makeatletter
 {\par\unskip\endMakeFramed%
 \at@end@of@kframe}
\makeatother

\definecolor{shadecolor}{rgb}{.97, .97, .97}
\definecolor{messagecolor}{rgb}{0, 0, 0}
\definecolor{warningcolor}{rgb}{1, 0, 1}
\definecolor{errorcolor}{rgb}{1, 0, 0}

\newcommand{\red}[1]{\textcolor[rgb]{0.90,0.00,0.10}{#1}}

\newcommand{\blue}[1]{\textcolor[rgb]{0.00,0.00,0.80}{#1}}

\newcommand{\Intv}[1]{\left[#1\right]}
\newtheorem{thm}{Theorem}

\newtheorem{defn}[thm]{Definition}

\newtheorem{remark}[thm]{Remark}

\newcommand{\ds}{\displaystyle}

\newcommand{\norm}[1]{\left\Vert#1\right\Vert}

\newcommand{\set}[1]{\left\{#1\right\}}

\newcommand{\rb}[1]{\left(#1\right)}

\newcommand{\R}{\mathbb{R}}

\newcommand{\Z}{\mathbb{Z}}

\newcommand{\Hp}{\mathbb{H}}

\newcommand{\N}{\mathbb{N}}

\numberwithin{equation}{section}

\usepackage{calrsfs}

\allsectionsfont{\normalsize\bfseries}

 \newcommand{\MC}[1]{ \mathcal{#1}}

\def\F{\mathbb F}

\numberwithin{equation}{section}
\allowdisplaybreaks

\usepackage[table, dvipsnames]{xcolor}
\DeclareMathAlphabet{\pazocal}{OMS}{zplm}{m}{n}
    \newcounter{example}[section]
\newenvironment{example}[1][]{\refstepcounter{example}\par\medskip
   \noindent \textbf{Example~\theexample. #1} \rmfamily}{\medskip}
\begin{document}
\title{
Topological comparison of  some dimension reduction methods  using persistent homology on EEG data
}
\author{Eddy Kwessi\footnote{Corresponding author: Department of Mathematics, Trinity University, 1 Trinity Place, San Antonio, TX 78212, Email: ekwessi@trinity.edu}}
\date{}
\maketitle
\begin{abstract}
In this paper, we explore how to use topological tools to compare dimension reduction methods. 
We first  make a brief overview of some of the methods often used  dimension reduction  such as Isometric Feature Mapping, Laplacian Eigenmaps, Fast Independent Component Analysis, Kernel Ridge Regression, t-distributed  Stochastic Neighbor Embedding. We then give a brief overview of some topological notions used in topological data analysis, such as, barcodes,  persistent homology, and Wasserstein distance.  Theoretically, these methods applied on a data set can be interpreted differently. From EEG data embedded into a manifold of high dimension, we apply these methods and we compare them across persistent homologies of dimension 0, 1, and 2, that is, across  connected components, tunnels and holes, shells around voids or cavities. We find that from three dimension clouds of points, it is not clear how distinct  from each other the methods are, but Wasserstein and Bottleneck distances, topological  tests of hypothesis,  and various methods show that the methods qualitatively and significantly differ across  homologies. 
\end{abstract}
\section{Introduction}
In topological data analysis, one is interested in understanding high dimensional structures from low dimensional ones and how discrete structures can be aggregated  to form a global structure. It can be a difficult task to even think or believe that high dimensional object exist beyond three dimensions since we can not visualize objects beyond a three-dimensional space. However, embedding theorems, for instance \cite{Whitney} and \cite{Takens1981} embedding theorems, clearly show that these high dimension structures do in fact exist. From a practical point of view, to make inferences on structures embedded in high dimensional ambient spaces,  some kind of dimensional reduction needs to occur. From a data analysis point of view, dimension reduction amounts to data compression where a certain amount of information may be lost. This dimension reduction is part of manifold learning, which can be understood as a collection of algorithms for recovering low dimension manifolds embedded into high dimensional ambient spaces, while preserving meaningful information, see \cite{MaFu}. The algorithms for dimension reduction may be classified into linear and nonlinear methods or parametric or nonparametric methods,  where the goal is to select or extract coarse features from high dimensional data. Among the pioneering linear methods is the Principal Component Analysis (PCA) introduced by \cite{Hotelling}. Its primary goal is to reduce the data to  a set of orthogonal linear projections ordered by decreasing  variances. Another linear method is the multidimensional scaling (MSD) where the data are aggregated using a measure of proximity, which could be a distance, or a measure of association such as correlation, or any other  method describing how close entities can be. Linear Discriminant Analysis (LDA) is a linear method similar to PCA consisting of writing a categorical dependent variable as a linear combination of continuous independent variables. As such, it is opposite to an Analysis of Variance (ANOVA) where the dependent variable is continuous and the independent variables are categorical. The focus of this paper will be on nonlinear techniques, which as their linear counterparts, aim to extract or select low dimensional features while preserving important information.  Since there are many such methods, our focus  will be on Isometric Feature Mapping (ISOMAP), Laplacian Eigenmaps, Fast Independent Component Analysis, (Fast-ICA), T-distributed Stochastic Neighbor Embedding (t-SNE). We will compare  them using  Persistent homology (PH). 
 PH is one the many techniques of topological data analysis (TDA) that can be used to identify features in data that remain persistent over multiple and different scales. This tool can provide new insights into seemingly known or unknown data and has the potential to uncover interesting  hidden information embedded within data. For instance, PH has been used to provide new insights on the  topology of deep neural networks, see \cite{Naizait2020}.  PH has successfully been used to provide new perspectives on viral evolution, see \cite{Chan2013}. The following examples of successful applications can be found in \cite{Otter2017} including but not  limited to  better understanding of sensor-network coverage, see \cite{Silva2007}, proteins, see \cite{Xia2014, Gameiro2015}, dimensional structure of DNA,  see \cite{Emmett2016}, see cells development, \cite{Rizvi2017}, robotics, see \cite{Bhattacharya2015,Pokorny2016, Vasudevan2013}, signal processing, see \cite{Chung2009, Guillemard2013}, spread of contagions, see \cite{Taylor2015}, financial networks, see \cite{Leibon2008}, applications in neuroscience, see \cite{Giusti2016,Sizemore2019}, time-series output of dynamical systems, see \cite{Maletic2016}, on EEG Epileply, see \cite{Chung2023}. The approach is the last reference is of particular interest to us. Indeed that paper,  authors considered EEG measured on the healthy person during sleep. They used the method of false nearest neighbors to estimate the embedding dimension. From there, persistent barcodes diagrams are obtained and reveal that topological noise persists at certain dimension and vanish some others. This paper has a similar approach  and is organized as follows: in Section \ref{sec2}, we review theories behind  some dimension reduction methods, then in Section \ref{sec3}, we give an overview of the essentials of persistent homology, in Section \ref{sec4}, we discuss how to apply persistent homology to data, and compare the methods on  an  EEG dataset, using persistent homology. Finally in Section \ref{sec5}, we make some concluding remarks. 
 \section{Review of a selected dimension reduction methods} \label{sec2}
 Let us note that some  of the review methods below are  extensively described in \cite{MaFu}. To have all of our ideas self-contained, let us re-introduced a few concepts.  In the sequel, $\norm{\cdot}$ is the euclidian norm in $\R^d$, for some $d\geq 3$. In the sequel, topological spaces $\MC{M}$ will considered to be second-countable Hausdorff, that is, (a) Every pair of distinct points has a corresponding pair of disjoint neighborhoods. (b) Its topology has a countable basis of open sets. This assumption is satisfied in most topological spaces and seems reasonable. 
 \subsection{Preliminaries}
 \begin{defn}
 A topological space   $\mathcal{M}$ is called a (topological) manifold if locally, it resembles a real $n$-dimensional Euclidian space, that is, there exists $n\in \N$ such that for all $x\in \mathcal{M}$, there exists a neighborhood $U_x$ of $x$  and a homeomorphism $f: U_x\to \R^n$. The pair $(U_x, f)$ is referred to as a chart on $\mathcal{M}$ and $f$ is called a parametrization at $x$.
 \end{defn}
 \begin{defn}
 Let $\mathcal{M}$ be a manifold. $\mathcal{M}$ is said to be smooth if given $x\in \mathcal{M}$, the parametrization  $f$ at $x$ has smooth or continuous partial derivatives of any order  and can be extended to a smooth function $F:\mathcal{M}\to \R^n$ such that $F\big|_{\mathcal{M}\cap U_x}=f$. 
 \end{defn}
 \begin{defn}  Let  $ \mathcal{M}$ and  $\mathcal{N}$ be  differentiable manifolds  and consider $\psi: \mathcal{M}\to \mathcal{N}$ a function. $\psi$ is said to be an  immersion if $\psi$ is a differentiable function and its  derivative is everywhere injective.
 In other words, $\psi:  \mathcal{M} \to \mathcal{N}$  is an immersion if $D_x\psi: T_x\MC{M}\to T_{\psi(x)}\MC{N}$ is an injective function at every point $x$ of $\MC{M}$, where $T_x \MC{M}$ represents the tangent plane to $\MC{M}$ at $x$.
 \begin{defn}
  \noindent Let  $ \mathcal{M}$ and  $\mathcal{N}$ be  differentiable manifolds. A function  $\psi: \mathcal{M}\to \mathcal{N}$ is an embedding  if $\psi$ is an injective immersion.
  \end{defn}


 \end{defn}
Let us introduce the notion of boundary of topological manifold that will be important in the sequel.
\begin{defn} Consider a Hausdorff topological manifold $\MC{M}$ homeomorphic to an open subsets of the half-euclidian space $\R^n_+$. Let the interior $Int(\MC{M})$ of $\MC{M}$ be the subspace of $\MC{M}$ formed by all points $s$ that have a neighborhood homeomorphic to $\R^n$. Then the boundary of $\MC{M}$ is defined as complement of $Int(\MC{M})$ in $\MC{M}$, that is, $\MC{M}\setminus Int(\MC{M})$, which is an $n-1$-dimensional topological manifold.
\end{defn}
  \subsection{Isomap} Isometric Feature Mapping (Isomap) was introduced by \cite{Tenenbaum2000}.  The data are considered to be a finite sample $\set{\bm{v}_i}$ from a smooth manifold $\mathcal{M}$. The two key assumptions are:  (a) there exists of an isometric embedding $\psi: \mathcal{M}\to \MC{X}$ where $\mathcal{X}=\R^d$, where the distance on $\MC{M}$ is the geodesic distance or the shortest curve connecting two points, (b) the smooth manifold $\MC{M}$ is a convex region of $\R^m$, where $m<<d$. The implementation phase has three main steps.
  \begin{enumerate}
  \item For a fixed integer $K$ and real number $\epsilon>0$, perform an $\epsilon-K$-nearest neighbor search using the fact that the geodesic distance $D^{\MC{M}}(v_i,v_j)$ between two points on $\MC{M}$ is the same (by isometry) as their euclidian distance $\norm{v_i-v_j}$ in $\R^d$.  $K$ is the number of data points selected within a ball of radius $\epsilon$.
  \item Having calculated the distance between points as above, the entire data set can be considered as a weighted graph with vertices $\bm{v}=\set{v_i}$ and edges $\bm{e}=\set{e_{ij}}$, where $e_{ij}$ connects $v_i$ with $v_j$ with a distance $w_{ij}=D^{\MC{M}}(v_i,v_j)$ considered as an  associated  weight. The geodesic distance between two data points $v_i$ and $v_j$ is estimated as the graph distance between the two edges, that is, the number of edges in a shortest path connecting them. We observe that this shortest path is found by  minimizing  the sum of the weights of its constituent edges. 
  \item Having calculated the geodesic distances $D^G=\set{w_{ij}}$ as above, we observe that $D^G$ is a symmetric matrix, so we can apply the classical Multidimensional Scaling algorithm (MDS) (see \cite{Torgerson1952}) to $D^G$ by mapping (embedding) them into a feature space $\MC{Y}$ of dimension $d$ while preserving the geodesic distance on $\MC{M}$. $\MC{Y}$ is generated by a $d\times m$ matrix whose $i$-th column represents the coordinates of $v_i$ in $\MC{Y}$.
  \end{enumerate}
  \subsection{Laplacian Eigenmaps}
  The Laplacian Eigenmaps (LEIM) algorithm was introduced by \cite{Belkin2002a}.  As above, the data ${\bm v}=\set{v_i}$ are supposed to be from a smooth manifold $\MC{M}$. It also has three main steps:
  \begin{enumerate}
  \item For a fixed integer $K$ and real number $\epsilon>0$, perform an $\epsilon-K$-nearest neighbor search on symmetric neighborhoods. Note that given two points $v_i, v_j$, their respective $K$-neighborhood $N_i^K$ and  $N_j^K$ are symmetric if and only $v_i\in N_j^K\Longleftrightarrow v_j\in N_i^K$. 
  \item For a given real number  $\sigma>0$ and  each pair of points $(v_i, v_j)$, calculate the weight $w_{ij}=e^{-\frac{\norm{u_i-v_j}^2}{2\sigma^2}}$ if $v_i\in N_j^K$ and $w_{ij}=0$ if $v_i\notin N_j^K$. Obtain the adjacency matrix $\bm{W}=(w_{ij})$. The data now form a weighted graph with vertices $\bm{v}$, with edges $\bm{e}=\set{e_{ij}}$,  and weights $\bm{W}=\set{w_{ij}}$, where $e_{ij}$ connects $v_i$ with $v_j$ with distance $w_{ij}$.
  \item Consider $\bm{\Lambda}=\set{\lambda_{ij}}$ be a diagonal matrix with $\ds \lambda_{ii}=\sum_{j} w_{ij}$ and define the graph Laplacian as  $\bm{L}=\bm{\Lambda}-\bm{W}$. Then $\bm{L}$ is positive definite so let $\widehat{\bm{Y}}$ be the $d\times n$ matrix  that minimizes $\ds \sum_{i,j}w_{ij}\norm{\bm{y}_i-\bm{y}_j}^2=tr(\bm{TLY^T})$. Then $\widehat{\bm{Y}}$ can used to embed $\MC{M}$ into a $d$-dimensional space $\MC{Y}$, whose $i$-th column represents the coordinates of $v_i$ in $\MC{Y}$.
  \end{enumerate}
   \subsection{Fast ICA}
 The Fast Independent Component Analysis (Fast-ICA) algorithms were introduced by \cite{Hyvarinen1999}. As above, the data $\bm{v}$ is considered to be from a smooth manifold $\MC{M}$. It is assumed that the data $\bm{v}$ is represented as an $n\times m$ matrix $(v_{ij})$ that can be flattened into a $ n\times m$ vector. As in Principal Component Analysis (PCA), in Factor Analysis, Projection Pursuit, or Independent Component Analysis (ICA), by considering the data as an $n\times m $-dimensional observed random variable, the goal  is  to determine a matrix $\bf{W}$ such that $\bm{s=W^Tv}$, where $\bm{s}$ is a $n\times m$-dimensional random variables having desirable properties such as optimal dimension reduction, or other interesting statistical properties such as minimal variance. Optimally, the components of $\bm{s}$ should provide source separation (the original data source $\bm{v}$ is assumed corrupted with noise) and feature extraction and be independent of each other. In a regular ICA, the matrix $\bm{W}$ is found by minimizing the mutual information, a measure of dependence between given random variables. In fast ICA algorithms, the matrix $\bm{W}$ is found by  using a Newton fixed point approach, with an objective function taken as   the differential entropy given as  $J_G(\bm{W})=\rb{\mathbb{E}[G(\bm{W^TW})]-\mathbb{E}[G(z)]}^2$, where it is assumed that $\bm{W}$ is such that $\mathbb{E}[\bm{(W^TW)}^2]=1$, and $z$ is standard normal distribution. $G$ is a function referred to as the  contrast function that include but is not limited to  $G(u)=\alpha^{-1}\log(\cosh(\alpha u)), G(u)=-\sigma^{-1}e^{-0.5\sigma u^2}, G(u)=0.25u^4$, where $\alpha\in [1,2]$ and $\sigma\approx 1$. From a dynamical system point of view, the fixed point is locally asymptotically stable with the exception of  $G(u)=0.25u^4$ where stability becomes global. For simplification purposes, let $g(x)=G'(x)$.  The key steps are:
 \begin{enumerate}
 \item Data preparation: it consists of centering the data $\bm{v}$ with respect to the column to obtain $\bm{v}^c$. That is, $\ds v^c_{ij}= v_{ij}-\frac{1}{m}\sum_{j=1}^mv_{ij}$, for $i=1,2,\cdots, n$. The  centered data are then  whitened, that is, $\bm{v}^c$ is linearly transformed  into $\bm{v}_w^c$,  a matrix of uncorrelated components. This is accomplished through  an eigenvalue decomposition of the covariance matrix $\bm{C=v^c(v^c)^T}$ to obtain two matrices $\bm{V}, \bm{E}$,  respectively of eigenvectors and eigenvalues so that $\bm{\mathbb{E}[C]=VEV^T}$.  The whitened data are found as $\bm{v^c_w=E^{-1/2}V^Tv^c}$ and simply referred to again as $\bm{v}$ for simplicity.
 \item Component extraction:  Let $F(\bm{W})=\mathbb{E}[\bm{v}g(\bm{W^Tv})]-\beta \bm{W}$ for a given constant $\beta=\mathbb{E}[\bm{W_a^T v}g(\bm{W_a^Tv})]$, where ${W_a}$ is the optimal weight matrix. Applying the Newton scheme ($x_{n+1}=x_n-F(x_n)[F'(x_n)]^{-1}$) to the  differentiable function $J_G$, we have 
 \begin{itemize}
 \item Select a random starting vector $\bm{W}_0$.
 \item For $n\geq 0$, $\bm{W}_{n+1}=\mathbb{E}[\bm{v}g(\bm{W}_n^T\bm{v})]-\mathbb{E}[g'(\bm{W}_n^T\bm{v})]\bm{W}_n$.
 \item   We then normalize $\bm{W}_{n+1}$ as $\ds \frac{\bm{W}_{n+1}}{\norm{\bm{W}_{n+1}}}$. 
 \item We repeat until a suitable convergence level is reached. 
 \item From the last matrix $\bm{W}$ obtained, we let $\bm{s=W^Tv}$.
 \end{itemize}
 \end{enumerate}
 \subsection{Kernel Ridge Regression}
The  Kernel Ridge Regression (KRR) is constructed as follows: as above, the data $\bm{v}$ is considered to be from a smooth manifold $\MC{M}$ of dimension say $d$. It is assumed that the data $\bm{v}$ is represented as an $n\times m$ matrix $\set{v_{ij}}$ that can be flattened into a $ n\times m$ vector. Suppose we are in possession of $\bm{u}=(u_1,u_2,\cdots, u_n)$ data corresponding to a response variable and covariates given as $\bm{v}=(\bm{v_1,v_2,\cdots, v_n})$ where $\bm{v_i}=(v_{ij})^T$ for $j=1,2,\cdots,m$.
 With the Least Square method, on can find the best linear model between the covariates $\bm{v}=(v_i)$ and the response $\bm{u}=(u_i)$, by minimizing the objective function $\ds L(\bm{W})=\frac{1}{2}\sum_{i=1}^L(u_i-\bm{W^Tv}_i)^2$, where $\bm{W}$ is a $1\times n$ vector. Least square methods are notorious for overfitting. The Ridge regression is a compromise that uses a penalized objective function such as $\ds L(\bm{W})=\frac{1}{2}\sum_{i=1}^L(u_i-\bm{W^Tv}_i)^2+\frac{\lambda}{2}\norm{\bm{W}}^2$. The solution can be found as $\ds \bm{W}=\rb{\lambda I+\sum_{i=1}^n\bm{v_iv_i}^T}^{-1}\rb{\sum_{i=1}^nu_i\bm{v}_i}$.
 In case the true nature of the relationship between the response and covariates is nonlinear, we can replace $\bm{v}_i$ with $\varphi(\bm{v}_i)$ where $\varphi$ is a nonlinear function $\R^m\to \R$. In particular, if the response is qualitative, that is, say labels, then we have a classification problem and $\varphi$ is referred to as feature map. Note that when using $\varphi$, the number of dimensions of  the problem is considerably high. Put  $\Phi=\varphi(\bm{v})=(\varphi(v_1),\varphi(v_2),\cdots, \varphi(v_n))$.  Replacing $\bm{v_i}$ with $\varphi(\bm{v_i})$, the solution above  becomes   $\ds \bm{W}=\rb{\lambda I+\sum_{i=1}^n\varphi(\bm{v_i})\varphi(\bm{v_i})^T}^{-1}\rb{\sum_{i=1}^nu_i\varphi(\bm{v}_i)}=(\lambda I+\Phi\Phi^T)^{-1}\Phi \bm{u}^T$. Consider  the following identity $AB^T(C+BAB^T)^{-1}=(A^{-1}+B^TC^{-1}B)^{-1}B^TC^{-1}$ for given invertible matrices $A,C$ and a matrix $B$. Applying this with $A=C=I$ and $B=\Phi$ we have   $\bm{W^T}=\bm{u}\Intv{\Phi^T(\lambda I+\Phi^T\Phi)^{-1}}=\bm{u}\Intv{(\lambda I+\Phi^T\Phi)^{-1}\Phi^T}$.
 Therefore, given a new value $\bm{v}_n$, the predicted value is $y_n=\bm{W}^T\Phi(\bm{v}_n)=\bm{u}(\Phi^T\Phi+\lambda I)^{-1}\Phi^T\Phi(\bm{v}_n)=\bm{u}(K+\lambda I)^{-1}\kappa(\bm{v}_n)$, where $\ds K=K(\bm{v}_i,\bm{v}_i)=\Phi^T\Phi=\sum_{i=1}^n\varphi(\bm{v}_i)^T\varphi(\bm{v}_i)$ and $\kappa(\bm{v}_n)=K(\bm{v}_i,\bm{v}_n)$. K is referred to as the kernel, which is the only quantity needed to be calculated, thereby significantly reducing the computational time and dimensionality of the problem. In practice, we may use a linear kernel $K(\bm{x},\bm{y})=\bm{x}^T\bm{y}$, or a Gaussian kernel $K(\bm{x},\bm{y})=e^{-\sigma\norm{\bm{x}-\bm{y}}^2}$, for some $\sigma>0$, where $\norm{\cdot}$ is a norm in $\R^m$  and $\sigma$ is given real constant.

 \subsection{t-SNE}
 Stochastic Neighbor Embedding (SNE) was proposed by \cite{Hinton2002}. t-SNE latter followed and was proposed by \cite{Maaten2008}.
 t-distributed stochastic neighbor embedding (t-SNE) is a dimension reduction method that amounts to assigning data two or three dimensional maps . As above, we consider the data $\bm{v}=(v_{ij})=(v_k)$ ($k=1,2,\cdots, N$ with $N=n\times m$) to be from a smooth  manifold $\MC{M}$ of high  dimension,  say $d$. The main steps of the method are:
 \begin{itemize}
 \item Calculate the asymmetrical probabilities $p_{kl}$ as $p_{kl}=\frac{e^{-\delta_{kl}}}{\sum_{k\neq l}e^{-\delta_{kl}}}$, where $\delta_{kl}=\frac{\norm{v_k-v_l}^2}{2\sigma_i}$ represents the dissimilarity between $v_k$ and $v_l$ and $\sigma_i$ is a parameter selected by the experimenter or by binary search. $p_{kl}$ represents the conditional probability that datapoint  $v_l$ is the neighborhood of datapoint $v_k$, if neighbors were selected proportionally to their probability density under a normal distribution centered at $v_k$ and variance $\sigma_i$.
 \item Assuming that the low dimensional data are $\bm{u}=(u_k),~k=1,2,\cdots, N$, the corresponding dissimilarity probabilities $q_{kl}$ are calculated under constant variance as $q_{kl}=\frac{e^{-d_{kl}}}{\sum_{k\neq l}e^{-d_{kl}}}$, where $d_{kl}=\norm{u_k-u_l}^2$ in the case of SNE and as $\ds q_{kl}=\frac{(1+d_{kl})^{-1}}{\sum_{k\neq l}{(1+d_{kl})^{-1}}}$, for t-SNE. 
 \item Then, we  minimize the Kullback-Leibler divergence between $p_{kl}$ and $q_{kl}$ given as $\ds L=\sum_{k=1}^N\sum_{l=1}^Np_{kl}\log\rb{\frac{p_{kl}}{q_{kl}}}$, using the gradient descent method with a momentum term with  the scheme $\ds \bm{w}^{t}=\bm{w}^{t-1}+\eta \frac{\partial L}{\partial \bm{u}}+\alpha(t)(\bm{w}^{t-1}-\bm{w}^{t-2})$ for $t=2,3,\cdots, T$ for some given $T$. Note that  $\bm{w}^{0}=(u_1,u_2,\cdots, u_{N})\sim N(0, 10^{-4}\bm{I})$, where $\bm{I}$ is the $N\times N$ identity matrix,  $\eta$ is a constant representing a learning rate, and $\alpha(t)$ is $t$-th  momentum iteration. We note that $\ds \frac{\partial L}{\partial \bm{u}}=\rb{\frac{\partial L}{\partial u_k}}$ for $k=1,2,\cdots, N$~ where $\ds \frac{\partial L}{\partial u_k}=4\sum_{l=1}^N(p_{kl}-q_{kl})(u_k-u_l)(1+d_{kl})^{-1}.$
 \item Then we use $\bm{u}=\bm{w}^T$ as the low dimensional representation of $\bm{v}$.
 \end{itemize}
\section{Persistent Homology} \label{sec3}
 In the sequel, we will introduce the essential ingredients needed to understand and compute persistent homology.
\subsection{Simplex complex}

\begin{defn}
A real $d$-simplex $S$ is a topological manifold of dimension $d$ that represents the convex hull of $d+1$ points. In other words:
\begin{equation}
S=\set{(t_0,t_1,\cdots, t_d)\in \R^d: t_i\geq 0 \quad \mbox{and $\ds \sum_{i=1}^dt_i=1$}}\;.
\end{equation}
\end{defn}
\begin{example} A $0$-simplex is a point, a $1$-simplex is an edge, a $2$-simplex is a triangle, a $3$-simplex is a tetahedron, a $4$-simplex is a pentachoron, etc.
\begin{figure}[H]
\centering\includegraphics[scale=.5]{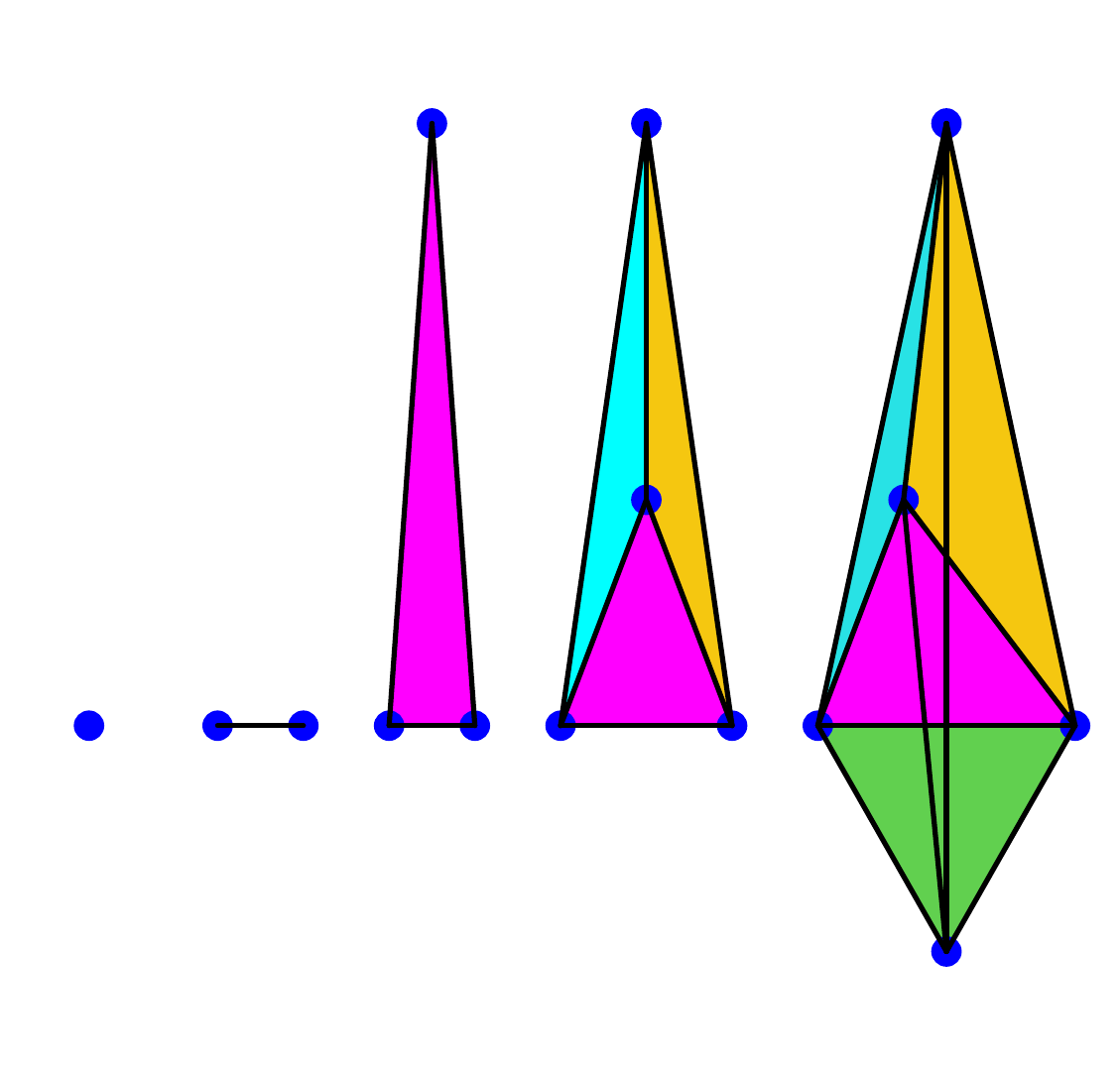}
\caption{An illustration of 0, 1, 2, 3, and 4-simplices.}
\end{figure}

\end{example}
\begin{remark}
We observe that a $d$-simplex $S$ can also be denoted as 
\[S=[V_0,V_1,\cdots, V_d],\quad \mbox{where $V_i=\set{\mbox{vertices of $V_i$ }}, i=0,1,\cdots, d$}\;.\]\\
We also note that the dimension of $V_i$ is $i$.
\end{remark}
\begin{defn} Given a simplex $S$, a face of $S$ is another simplex $R$ such that $R\subseteq S$ and such that the vertices of $R$ also the vertices of $S$.
\end{defn}
\begin{example} Given a 3-simplex (a tetrahedron), it has 4  different 2-simplex or 2 dimensional faces, each of them with three 1-simplex or 1-dimensional faces, each with three 0-simplex or 0-dimensional faces.  
\end{example}
\begin{defn} A simplicial complex $\Sigma$ is a topological space formed by different simplices not necessarily of the same dimension which have to satisfy the gluing condition, that is: 
\begin{enumerate}
\item Given $S_i\in \Sigma$, its face $R_i\in \Sigma$. 
\item Given $S_i, S_j\in \Sigma$, either $S_i\cap S_j=\emptyset$ or $S_i\cap S_j=R_i=R_j$, the  faces of $S_i$ and $S_j$ respectively.
\end{enumerate}
\end{defn}
\begin{figure}[H]
\centering \includegraphics[scale=0.5]{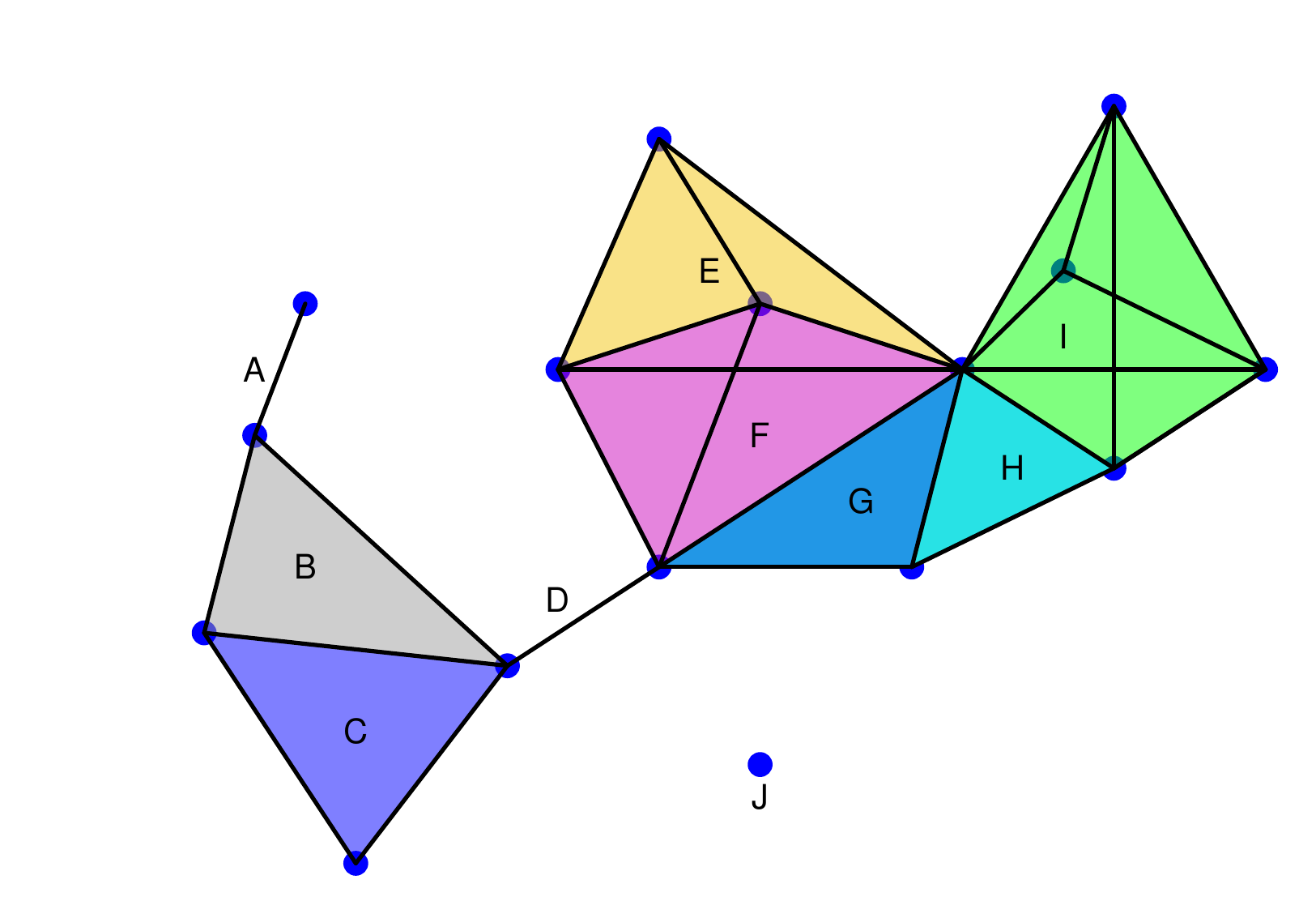}
\caption{Example of a simplicial complex. $J$ is 0-simplex,  $A$ and $D$ are 1-simplices, $B, C, G$, and $H$ are 2-simplices, $E$ and $F$ are 3-simplices, and $I$ is a 4-simplex. We note that $A\cap B$ is a 0-simplex. $B\cap C$ is a 1-simplex and a face of $B$ and $C$ respectively. $E\cap F$ is a 2-simplex and a face of $E$ and $F$. $G\cap H$ is a 1-simplex and $I \cap H$ is 1-simplex.}
\end{figure}
\noindent It is important to observe that a simplicial complex can be defined very abstractly. Indeed,
\begin{defn} A simplicial complex $\Sigma=\set{S: S\subseteq \Omega}$ is a collection of non-empty subsets of a set $\Omega$ such that 
\begin{enumerate}
\item For all $\omega\in \Omega$, then $\set{\omega}\in \Sigma$.
\item For any set $U$  such that $S\subset U$ for some $S\in \Sigma$, then $U\in \Sigma$.
\end{enumerate}
\end{defn}
\begin{example} Let $\Omega=\set{1,2,3,4}$. We can define the following simplicial complexes on $\Omega$.\\
1. $\ds \Sigma_1=\set{\set{1},\set{2},\set{3},\set{4}, \set{1,2},\set{1,3},\set{2,3},\set{1,2,3}}\;.$\\
2. $\Sigma_2=\MC{P}(\Omega)\setminus \set{\emptyset}$, where $\MC{P}(\Omega)$ is the set of all subsets of $\Omega$.
\end{example}

\subsection{Homology and persistent homology}
\begin{defn} Let $\Sigma$ be a simplicial complex.  \\
\noindent We define the Abelian group generated by the $j$-simplices of $\Sigma$ as $C_j(\Sigma)$.\\
\noindent We define a boundary operator associated with $C_j(\Sigma)$ as a homomorphism \[\partial_j: C_j(\Sigma)\to C_{j-1}(\Sigma)\;.\]
\noindent We define the chain complex associated with $\Sigma$ as the collection of pairs \[C(\Sigma)=\set{(C_j(\Sigma),\partial_j), j\in \Z}\;.\]
\end{defn}
\noindent Now we can define a homology group associated with a simplicial complex.
\begin{defn} Given a simplicial complex $\Sigma$, put $A_j(\Sigma):=kern(\partial_j)$ and $B_j(\Sigma):=Im(\partial_{j+1})$. Then  the $j$th homology group $\mathbb{H}_j(\Sigma)$of $\Sigma$ is defined as quotient group between $A_j(\Sigma)$ and $B_j(\Sigma)$, that is, 
\[\Hp_j(\Sigma)=\frac{A_j(\Sigma)}{B_j(\Sigma)}\;.\]
What this reveals is the presence of ``holes" in a given shape. 
\end{defn}
\begin{remark}
It is important to observe that $\ds \Hp_j(\Sigma)=\frac{<\mbox{$j$-dimensional cycles}>}{<\mbox{$j$-dimensional boundaries}>}$, where $<U>$ stands for the span of $U$, and a cycle is simply a shape similar to a loop but without necessarily a starting point.\\
\noindent Another important remark is that the boundary operator can indeed be defined as 
\[\partial_j(\Sigma):=\sum_{k=0}^j(-1)^k[V_0,\cdots, \widehat{V}_{-k},\cdots, V_j]\;,\]
where $\widehat{V}_{-k}$ means not counting the vertices of $V_k$. This  shows that $\partial_j(\Sigma)$ lies in a $j-1$-simplex. \\
Another remark is that $\partial_{j-1}\circ\partial_j=0$ for $0\leq j\leq d$\;.
\end{remark}
Now that we know that homology reveals the presence of ``holes", we need to find a way of determining how to count these ``holes".
\begin{defn} Given a simplicial complex $\Sigma$, the $j$th Betti number $b_j(\Sigma)$ is the rank of $\Hp_j(\Sigma)$ or \[b_j(\Sigma)=dim(A_j(\Sigma))-dim(B_j(\Sigma))\;.\]
 In other words, the smallest cardinality of a generating set of  the group $\Hp_j(\Sigma)$. \\
 In fact since the elements of $A_j(\Sigma)$ are $j$-dimensional cycles and that of $B_j(\Sigma)$ are $j$-dimensional boundaries,  the Betti number counts the number of independent $j$-cycles not representing the boundary of any collection of simplices of $\Sigma$.
\end{defn}
\begin{example} 
\begin{enumerate}
\item $b_0$ is the number of connected components of the complex.
\item \noindent $b_1$ is the number of tunnels and holes.
\item \noindent $b_2$ is the number of shells around cavities or voids.
\end{enumerate}
\end{example}
\begin{defn} Let  $\Sigma$ be simplicial complex and  let $N$ be a positive integer. A filtration of $\Sigma$ is a a nested  family $\Sigma^F_N:=\set{\Sigma_p, 0\leq p\leq N}$ of sub-complexes of $\Sigma$ such that 
\[\Sigma_0\subseteq \Sigma_1\subseteq \Sigma_2\subseteq\cdots\subseteq \Sigma_N=\Sigma\;.\]
\end{defn}
Now let $\F_2$ be the field with two elements and let $0\leq p\leq q\leq N$ be two integers. Since $\Sigma_p\subseteq \Sigma_q$, the inclusion map $Incl_{pq}: \Sigma_p\to \Sigma_q$ induces an $\F_2$-linear map $g_{pq}: \Hp_j(\Sigma_p) \to \Hp_j(\Sigma_q)$. We can now define, for any $0\leq j\leq d$,  the $j$-th persistent homology of a simplicial complex $\Sigma$.
\begin{defn} Consider  a simplicial complex $\Sigma$ with filtration $\Sigma^F_N$, for some positive integer $N$. The $j$-th persistent homology $\Hp_j^{p\to q}(\Sigma)$ of $\Sigma$ is defined as the pair:
\[\Hp_j^{p\to q}(\Sigma,\F_2):=(\set{\Hp_j(\Sigma_p), 0\leq p\leq N}, \set{g_{pq}, 0\leq p\leq q\leq N})\;.\]
\end{defn}
\noindent In a sense, the $j$-th persistent homology provides  a more refined information than the homology of the simplicial complex in that it informs us of the changes of features such as connected components, tunnels and holes, shells around voids  through the filtration process. It can be visualized using a ``barcode" or a persistent diagram.  The following definition is borrowed from \cite{Ghrist2008}:
\begin{defn} Consider a simplicial complex $\Sigma$, a positive integer $N$, and  two integers $0\leq p\leq q\leq N$.  The barcode of the $j$-th persistent homology $\Hp_j^{p\to q}(\Sigma,\F_2)$ of $\Sigma$ is a graphical representation of $\Hp_j^{p\to q}(\Sigma,\F_2)$ as a collection of horizontal line segments in a plane whose horizontal axis corresponds to a parameter and whose vertical axis represents an arbitrary ordering of homology generators.
\end{defn}
\noindent We finish this section  with the introduction of the Wasserstein and Bottleneck distances, used for the comparison of  persistent diagrams. 
\begin{defn} Let $p>1$ be a real number. Given two persistent diagrams $X$ and $Y$, the $p$-th Wasserstein distance $W_p(X,Y)$ between $X$ and $Y$ is defined as 
\[W_p(X,Y):=\inf_{\eta: X\to Y}\sum_{x\in X}\norm{x-\eta(x)}_{\infty}^p\;,\]
where $\eta$ is a perfect matching between the intervals of $X$ and $Y$.\\
The Bottleneck distance is obtained when $p=\infty$, that is, it is given as 
\[W_{\infty}(X,Y):=\inf_{\eta: X\to Y}\sup_{x\in X}\norm{x-\eta(x)}_{\infty}\;.\]

\end{defn}
\section{Applications to Data }  \label{sec4}
In the presence of data, simplicial complexes will be replaced  by sets of data indexed by a parameter, therefore transforming these sets into parametrized topological entities. On these parametrized topological entities, the notions of persistent homology  introduced above can be computed, especially the Betti number, in the form of ``barcode". To see how this could be done, let us consider the following definitions:
\begin{defn} For a given collection of points $\set{s_{\delta}}$ in a manifold $\MC{M}$ of dimension $n$, its   \v{C}ech complex $C_{\delta}$ is a simplicial complex formed by $d$-simplices obtained from a sub-collection $\set{x_{\delta, k}, 0\leq k\leq d, 0\leq d\leq n}$ of points such that taken pairwise, their $\delta/2$-ball neighborhoods have a point in common.
\end{defn}
\begin{defn}For a given collection of points $\set{s_{\delta}}$ in a manifold $\MC{M}$ of dimension $n$, its Rips complex $R_{\delta}$ is a simplicial complex formed by $d$-simplices obtained from a sub-collection  $\set{x_{\delta, k}, 0\leq k\leq d, 0\leq d\leq n}$ of points which are pairwise within a  distance of $\delta$.
\end{defn}
\begin{remark} 
\noindent 1. It is worth noting that in practice, Rips complexes  are easier to compute than \v{C}ech complexes, because the exact definition of the distance on $\MC{M}$ may not be known. \\
\noindent 2. More importantly, from a data analysis point of view, Rips complexes are good approximations (estimators) of \v{C}ech complexes. Indeed, a result from \cite{Silva2007} shows that  given $\delta>0$, there exits a chain of  inclusions $R_{\delta}\hookrightarrow C_{\delta/\sqrt{2}}\hookrightarrow R_{\delta/\sqrt{2}}$.
\noindent 3. Though Rips complexes and barcodes seems like challenging objects to wrap one's head around, there is an ever growing list of algorithms from various languages that can be used  for their visualization. All the analysis below has been done using R, in particular the TDA package in R.
\end{remark}
\subsection{Randomly generated data}
We generated 100  data points  sampled randomly in the square $[-5,5]\times [-5,5]$. In Figure \ref{fig3} and Figure \ref{fig4} below, we illustrate the Rips and barcode changes through a filtration. 

\begin{figure}[H]
\centering \includegraphics[width=1\textwidth, height=.5\textwidth]{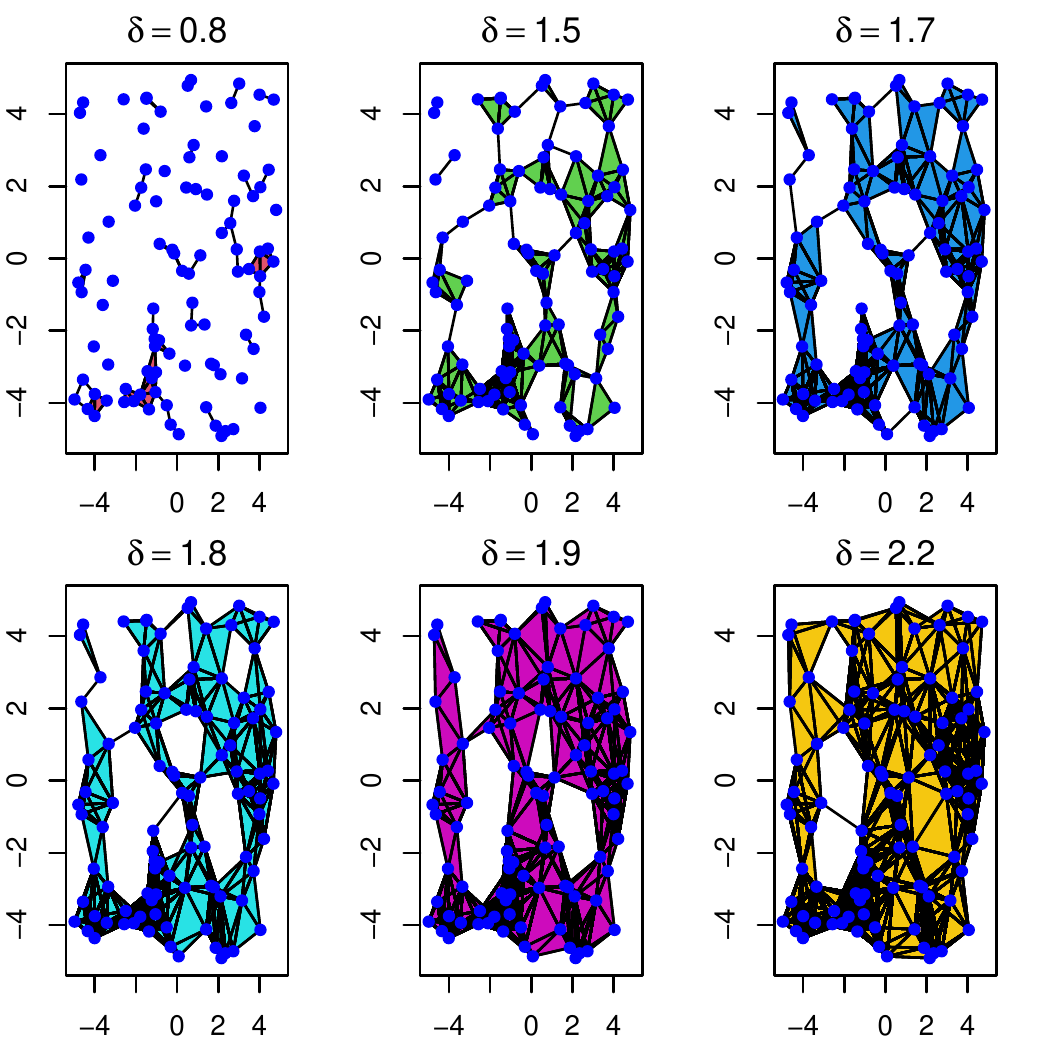}
\caption{Example of the evolution of Rips Complexes $\set{R_{\delta}}$ through a filtration with parameter $\delta$.  As we move from left to right, it shows how samples points (blue dots), first form 0-simplices, then 1-simplices, and so on. In particular, it shows how connected components  progressively evolve to form different types of holes.}
\label{fig3}
\end{figure}

\begin{figure}[H]
\centering \includegraphics[width=1\textwidth, height=.5\textwidth]{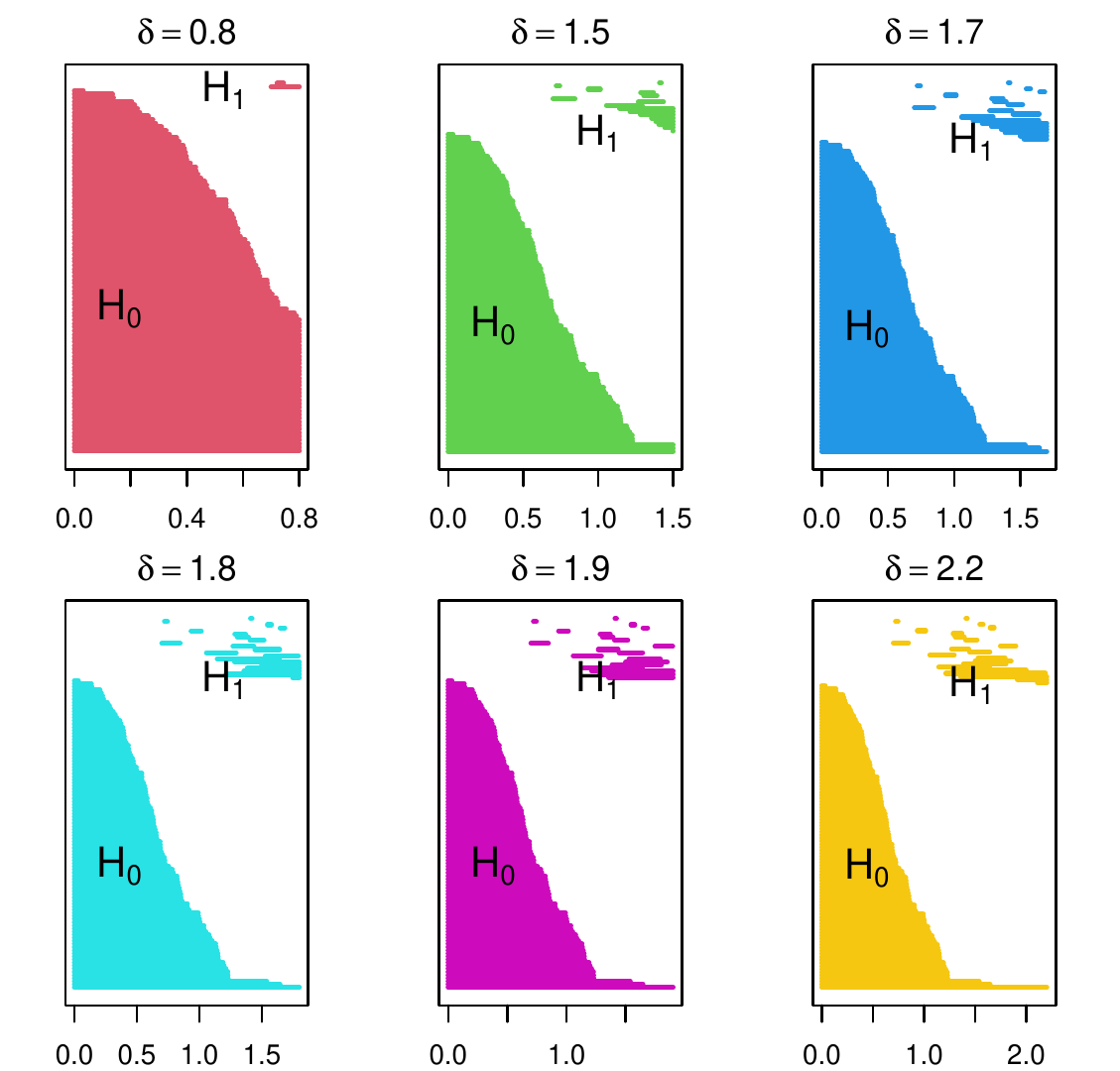}
\caption{Example of the evolution of barcodes  through a filtration with parameter $\delta$ for the same data as above. As we move from left to right, from top to bottom, it shows appearance and disappearance of lines ($\Hp_0$) and holes ($\Hp_1$) as the parameter $\delta$ changes.  It shows that certain lines and holes persist through the filtration process.}
\label{fig4}
\end{figure}
\subsection{EEG epilepsy data}
\subsubsection{Data Description}The main purpose of the manuscript it to analyze EEG data. We will consider a publicly available (at \url{http://www.meb.unibonn.de/epileptologie/science/physik/eegdata.html}) epilepsy data set called here EDATA for simplicity. The data consist of five sets  A, B, C, D, and E. Each containing 100 single-channel EEG segments of 23.6 seconds, each of which was selected  after visual inspection for artifacts (such as acoustic and electrical shielding, separation of earth ground for laboratory, interconnectivity of devices on the same phase and ground centrally and locally) and has passed a weak stationarity criterion. Sets A and B  were  obtained from surface EEG recordings of five healthy subjects  with eyes open and closed, respectively. Data  were obtained  in seizure-free intervals from five patients in the epileptogenic zone for set D and from the hippocampal formation of the opposite hemisphere of the brain for set C. Set E contains seizure activity, selected from all recording sites exhibiting ictal activity. Sets A and B have been recorded extracranially, whereas sets C, D, and E have been recorded intracranially. All EEG signals were recorded with the same 128-channel amplifier system, using an average common reference omitting electrodes containing pathological activity (C,D, and E) or strong eye movement artifacts (A and B). After 12 bit analog-to-digital conversion, the data were written continuously onto the disk of a data acquisition computer system at a sampling rate of 173.61 Hz. Band-pass filtersettings were 0.53--40 Hz (12 dB/oct.)
\subsubsection{Data analysis}
The approach is  to first  embed the data into a manifold of high dimension. This was already done in \cite{Kwessi2021}.  The dimension $d=12$ was found using the method of false nearest neighbors. Depending on the set used,  the size of the data can be very  large: for example ($4097\times 100\times 5=2,048,500$) making it very challenging to analyze holistically. In \cite{Kwessi2021}, we proposed to construct a complex structure (using all  100 channels for all 5 groups) whose volume changes per group. We would like to analyze the data further from a persistent homology point of view. 
This would mean analyzing 500 different persistent diagrams  and making  an  inference. We note that simplicial complexes of this data sets are very large  (2 Millions+). Fortunately, we can use the Wasserstein distance to compare persistent diagrams. To clarify, we will use each of the dimension reduction method introduced earlier, then proceed with construction of persistent diagrams. We will then compare them by method and by sets (A, B, C, D, and E).\\

\noindent \textbf{Single-channel Analysis:}\\

\begin{figure}[H]
\centering \includegraphics[width=1\textwidth, height=.6\textwidth]{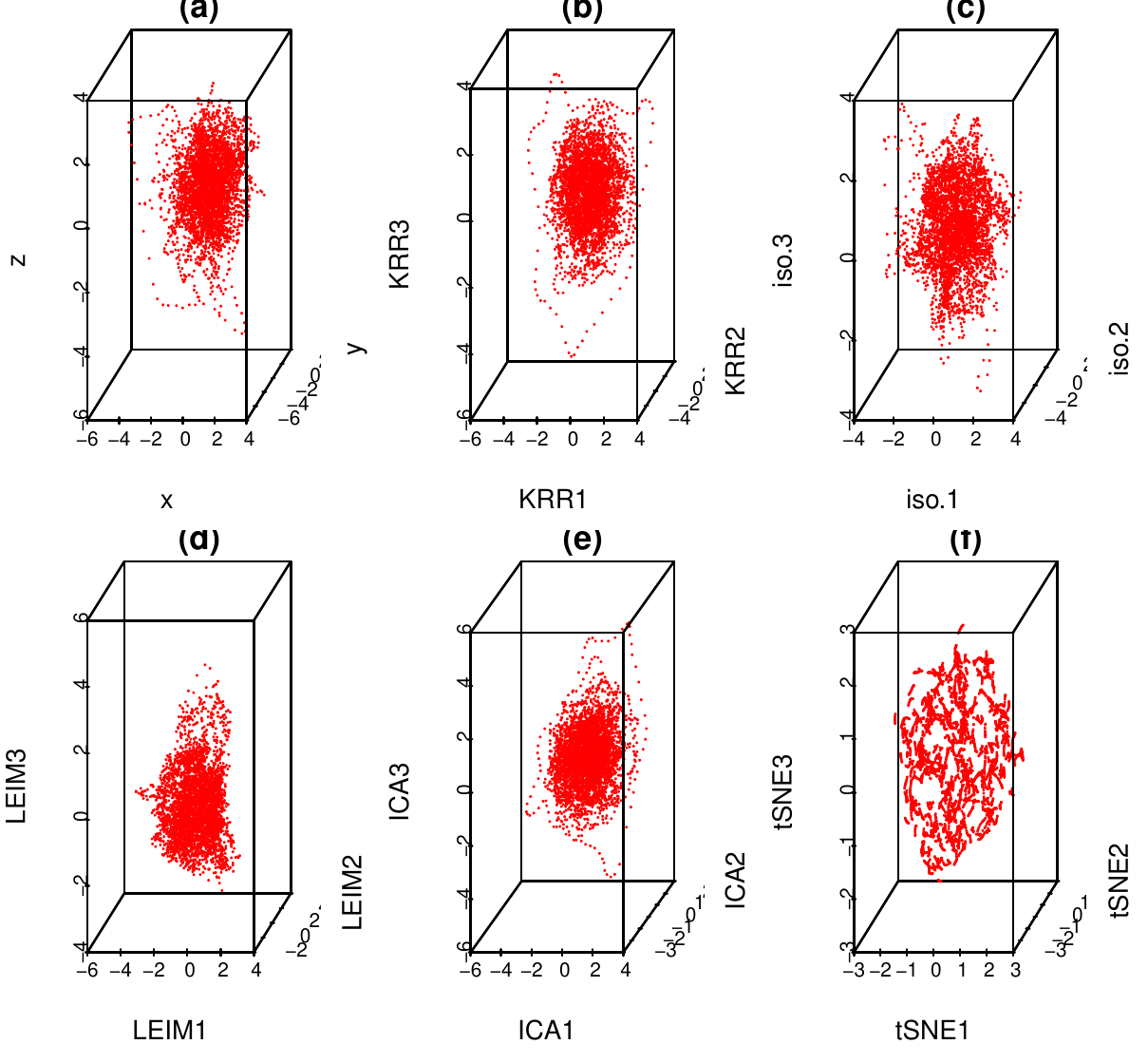} %
\caption{Scatterplots for a Takens projection method (a), KRR method (b), Isomap (c), LEIM (d), ICA (e), and t-SNE (f).}
\label{fig5}

\end{figure}

\begin{figure}[H]
\centering \includegraphics[width=1\textwidth, height=.6\textwidth]{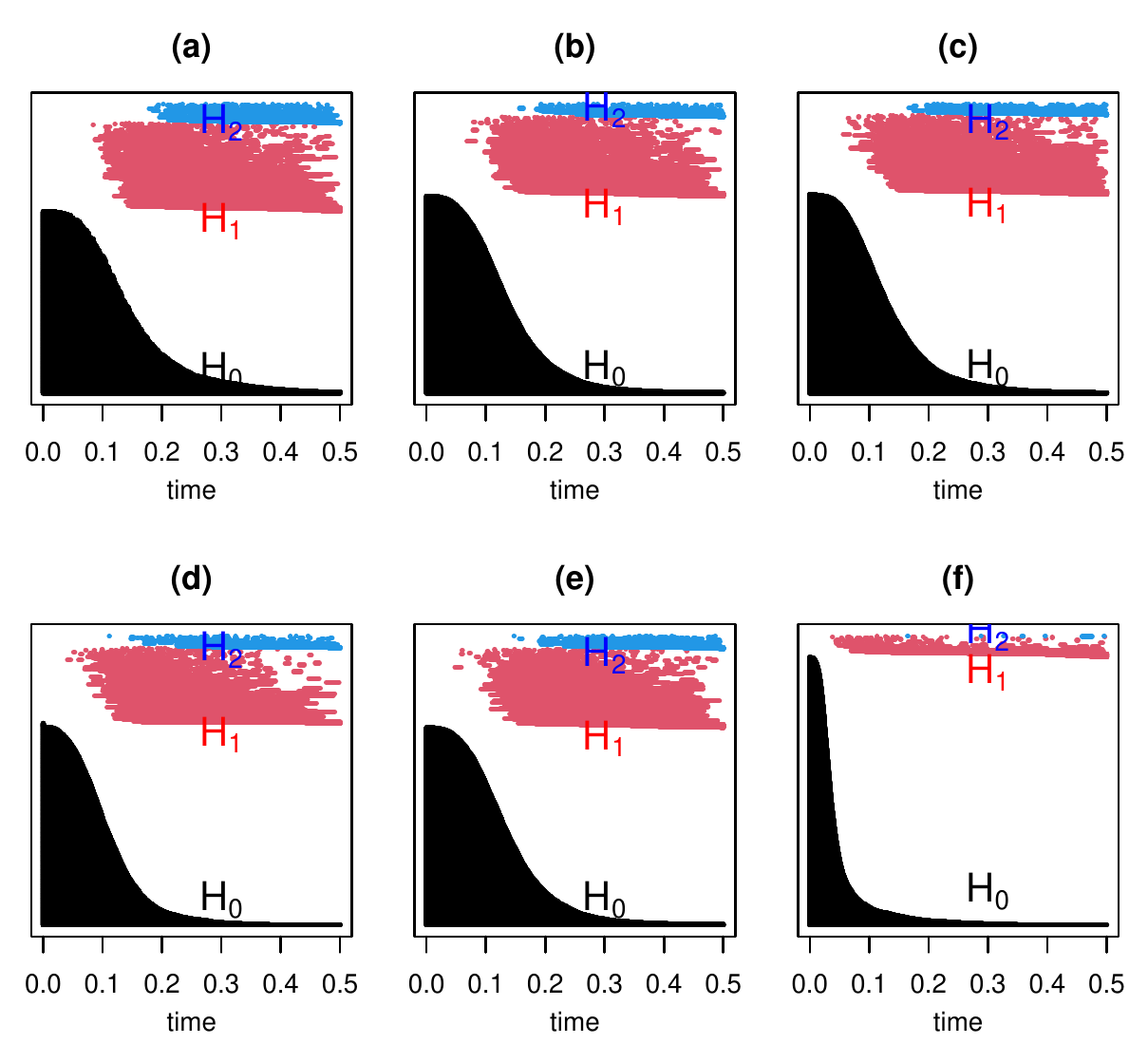}
\caption{Barcodes  for a Takens projection method(a), KRR method (b), Isomap (c), LEIM (d), ICA (e), and t-SNE (f). }
\label{fig6}
\end{figure}

\noindent  Suppose we select at random  one channel among the 100  from set D. Figure \ref{fig5} below represents a 3 dimensional representation of the embedded data using Takens embedding method (Tak), plotted using the 3 first three delayed coordinates $x=x(t), y=x(t-\rho), z=x(t-2\rho)$ where $\rho=1\Delta t$, with $\Delta t=\frac{1}{fs}=5.76$ ms in Figure \ref{fig5} (a), then the first three coordinates in the case of Kernel Ridge Regression (KRR$_i$) in Figure \ref{fig5} (b), Isomap (iso.$i$) in Figure \ref{fig5} (c), Laplacian Eigen Maps (LEIM$_i$) in Figure \ref{fig5} (d), Fast  Independent Component Analysis (ICA$_i$) in Figure \ref{fig5} (e), and t-Distributed Stochastic Neighbor Embedding (t-SNE$_i$) Figure \ref{fig5} (f). From these three dimensional scatter plots, we can visually observe that the  t-SNE plot (Figure \ref{fig5} (f)) is  relatively different from the other five since it seems to have more larger  voids. How different is difficult to tell from the naked-eye.
Figure \ref{fig6} represents their corresponding barcodes. It is much clearer looking at the the persistent diagram for t-SNE (Figure \ref{fig6} (f)) that it  is very different from the other five, when looking at $\Hp_0, \Hp_1$ and $\Hp_2$. Now, a visual comparison is not enough to really assert a significant difference. Using the Bottleneck distance, we calculate the distance between the respective persistent diagrams, for $\Hp_0$ and $\Hp_1$ in Table \ref{table1} (a)  and $\Hp_2$ in Table \ref{table1} (b)  below. We observe  from the first table that the Bottleneck distance at $\Hp_0$ and $\Hp_2$ for t-SNE are almost twice as large as for the other methods. They are comparable to that of LEIM at $\Hp_1$. The other methods have comparable Bottleneck distances at $\Hp_0, \Hp_1$, and $\Hp_2$, confirming what we already suspected visually in Figure \ref{fig5} and Figure \ref{fig6}. 


\begin{table}[H]
\resizebox{1\textwidth}{!}{\begin{minipage}{1.8\textwidth}
\centering
\begin{tabular}{cc}
{\bf (a)} & {\bf (b)}\\
\begin{tabular}{|l|r|r|r|r|r|r|}
\hline  
 $\Hp_0$ & Tak & Iso & KRR & ICA & LEIM & TSNE\\
\hline
Tak &  &  & & && \\
\hline
Iso & 0.0945019 &  & & & & \\
\hline
KRR & 0.0957546 & 0.0200035 & &  &   & \\
\hline
ICA & 0.0982795 & 0.0157002 & 0.0071899 & && \\
\hline
LEIM & 0.1678820 & 0.1182656 & 0.1247918 & 0.1205499 &  & \\
\hline
TSNE & 0.2238167 & 0.1730406 & 0.1817924 & 0.1759454 & 0.1162392 & \\
\hline
\end{tabular} &
\begin{tabular}{|l|r|r|r|r|r|r|}
\hline
  \red{$\Hp_1$}/\blue{$\Hp_2$}  & Tak & Iso & KRR & ICA & LEIM & TSNE\\
\hline
Tak &  & \blue{0.0363205} & \blue{0.0301992} & \blue{0.0292631} & \blue{0.0291247} & \blue{0.0551774}\\
\hline
Iso & \red{0.0340282} &  & \blue{0.0330687} & \blue{0.0290406} & \blue{0.0236890} & \blue{0.0598517}\\
\hline
KRR & \red{0.0317261} & \red{0.0279460} &  & \blue{0.0207599} & \blue{0.0212138} & \blue{0.0647935}\\
\hline
ICA & \red{0.0310771} & \red{0.0270919} & \red{0.0208086} &  & \blue{0.0242277} & \blue{0.0611090}\\
\hline
LEIM & \red{0.0607389} & \red{0.0725585} & \red{0.0702695} & \red{0.0682761} &  & \blue{0.0542615}\\
\hline
TSNE & \red{0.0757815} & \red{0.0959521} & \red{0.0864587} & \red{0.0861522} & \red{0.0785030} & \\
\hline
\end{tabular}
\end{tabular}
\end{minipage}}
\caption{Bottleneck distance between the persistent diagrams above  at $\Hp_0$ (black), at $\Hp_1$ (blue), and at $\Hp_2$ (red)\;.}
\label{table1}
\end{table}

\noindent The analysis above was done using  a single channel, selected at random from the set D. It seems to suggest that the t-SNE method is  different from the other five dimension reduction methods discussed above. Strictly speaking, non zero Bottleneck distances are indication of structural topological differences. What they do not say however is, if the differences observed  are significant. To address the issue of significance, we will perform a pairwise permutation test. Practically, from set $j$ and  channel $i$, we will obtain a persistent diagram $\MC{D}^{(j)}_i\sim \MC{P}^{(j)}$ where $j\in \set{1, 2, 3, 4, 5}, i\in \set{1,2,\cdots, 15}$, and $\MC{P}^{(j)}$ is the true underlying distribution of persistent diagrams, see \cite{Mileyko2011} for the  existence of these distributions. We will conduct a pairwise permutation test with null hypothesis  $H_0: \MC{P}^{(j)}=\MC{P}^{(j')}$ and alternative hypothesis $H_1: \MC{P}^{(j)}\neq \MC{P}^{(j')}$. We will use landscape functions (see \cite{Berry2020}) to obtain test statistics. The p.values obtained were found to be  very small, suggesting that the differences above  are indeed all significative across  $\Hp_0, \Hp_1$, and $\Hp_2$. \\

\noindent \textbf{Multiple-channel Analysis:}\\

\noindent (a) \noindent {Within set analysis}\\

\noindent In each set, we make a random selection of 15 channels, and we compare the Bottleneck distances obtained. This means having 15 tables of distance such as Table \ref{table1} (b) above. There will be consistency if the cell  value  $k (i,j)$ in  Table $k$, where $k\in \set{1,2,\cdots, 15}$ and $i,j\in \set{1,2,3,4,5}$ is barely different from $k'(i,j)$ of Table $k'$. Large differences will be an indication of topological differences between the methods within the sets. In Figure \ref{fig7} below, the $y$-axis represents Bottleneck distances and the $x$-axis represents channels indices. The red color is indicative of the Bottleneck distance between persistent diagrams on $\Hp_1$ and the blue color on $\Hp_2$ from data generated from each of the methods above. We see that overall, while there are small fluctuations from channels to channels  on $\Hp_1$, the largest fluctuations actually occur on $\Hp_2$. A deeper analysis reveals that in fact, the large fluctuations are due to large distance between t-SNE and the other five methods.
This confirms the earlier observations (refer to Figure  \ref{fig6} and Table \ref{table1} above) that persistent diagrams are really different on $\Hp_2$. Topologically, this means that shells around cavities or voids  that persist are not the same when using different dimension reduction methods. However,  the small fluctuations on $\Hp_1$ do not mean that tunnels and holes that persist are the same. Rather,  they do indicate is that they may not be all very different. 
 \\
\begin{figure}[H]
\resizebox{0.8\textwidth}{!}{\begin{minipage}{2\textwidth}
\centering
\begin{tabular}{ccc}
(A) &  (B) & (C)\\
\includegraphics{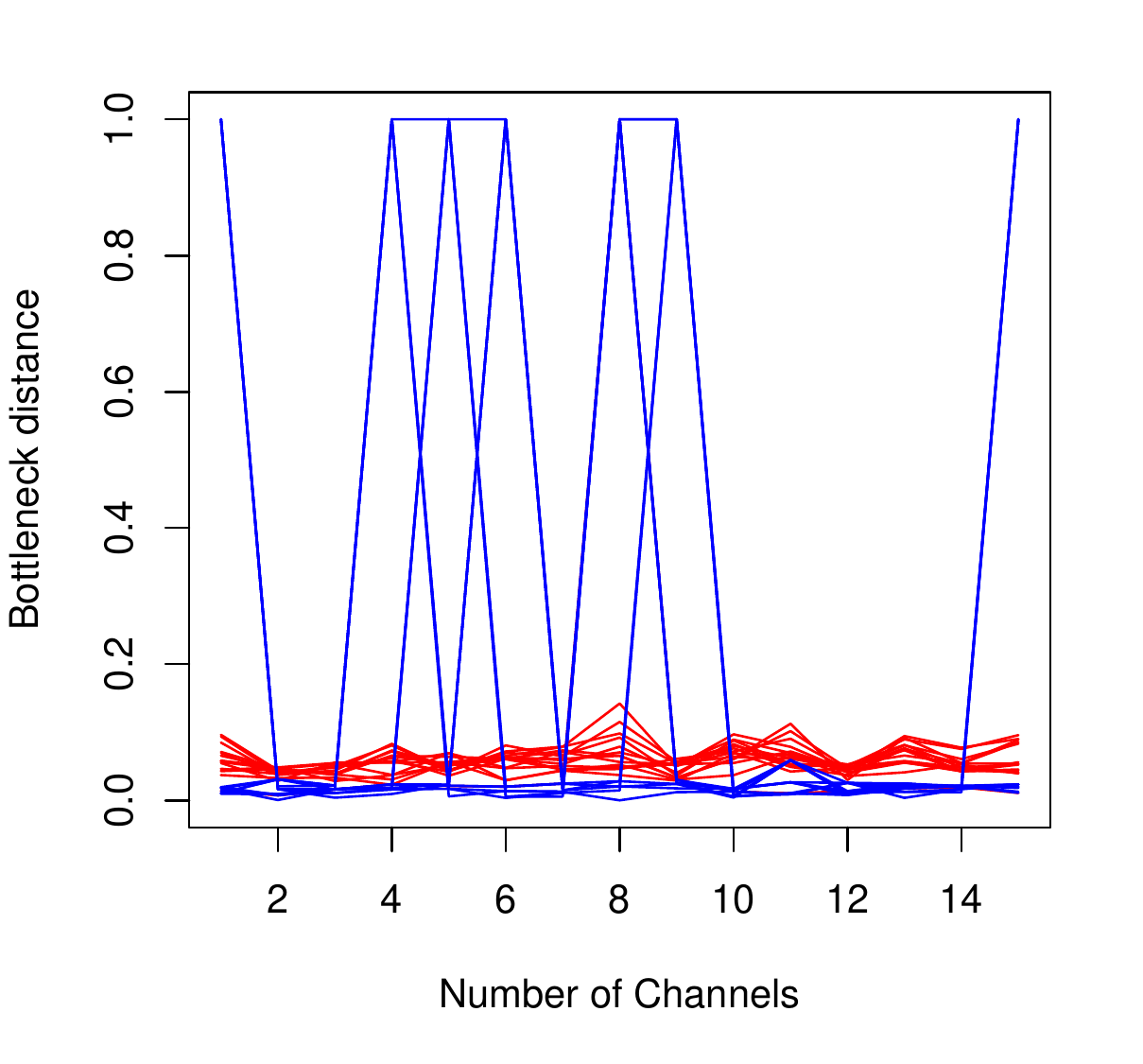} & \includegraphics{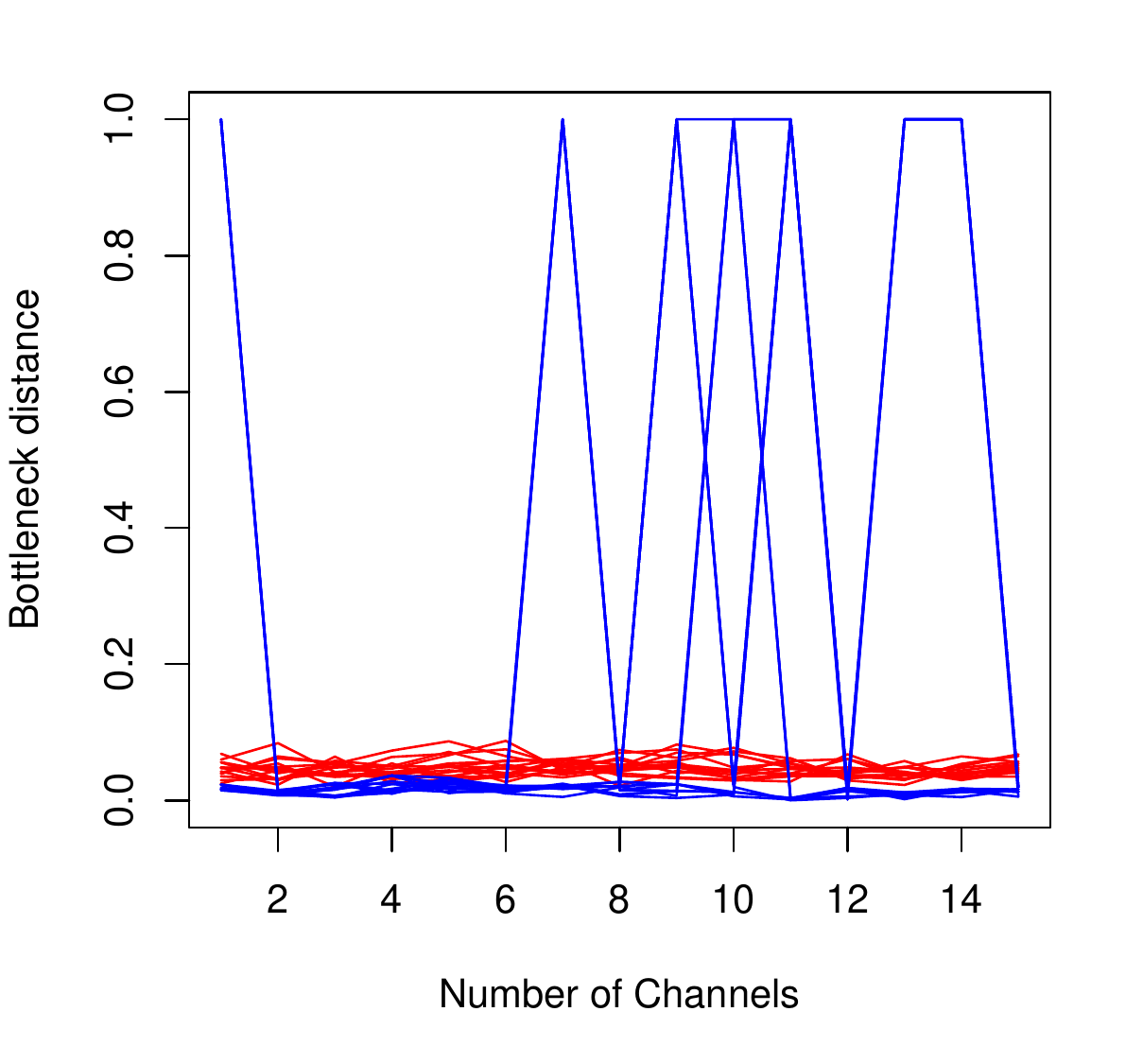} & \includegraphics{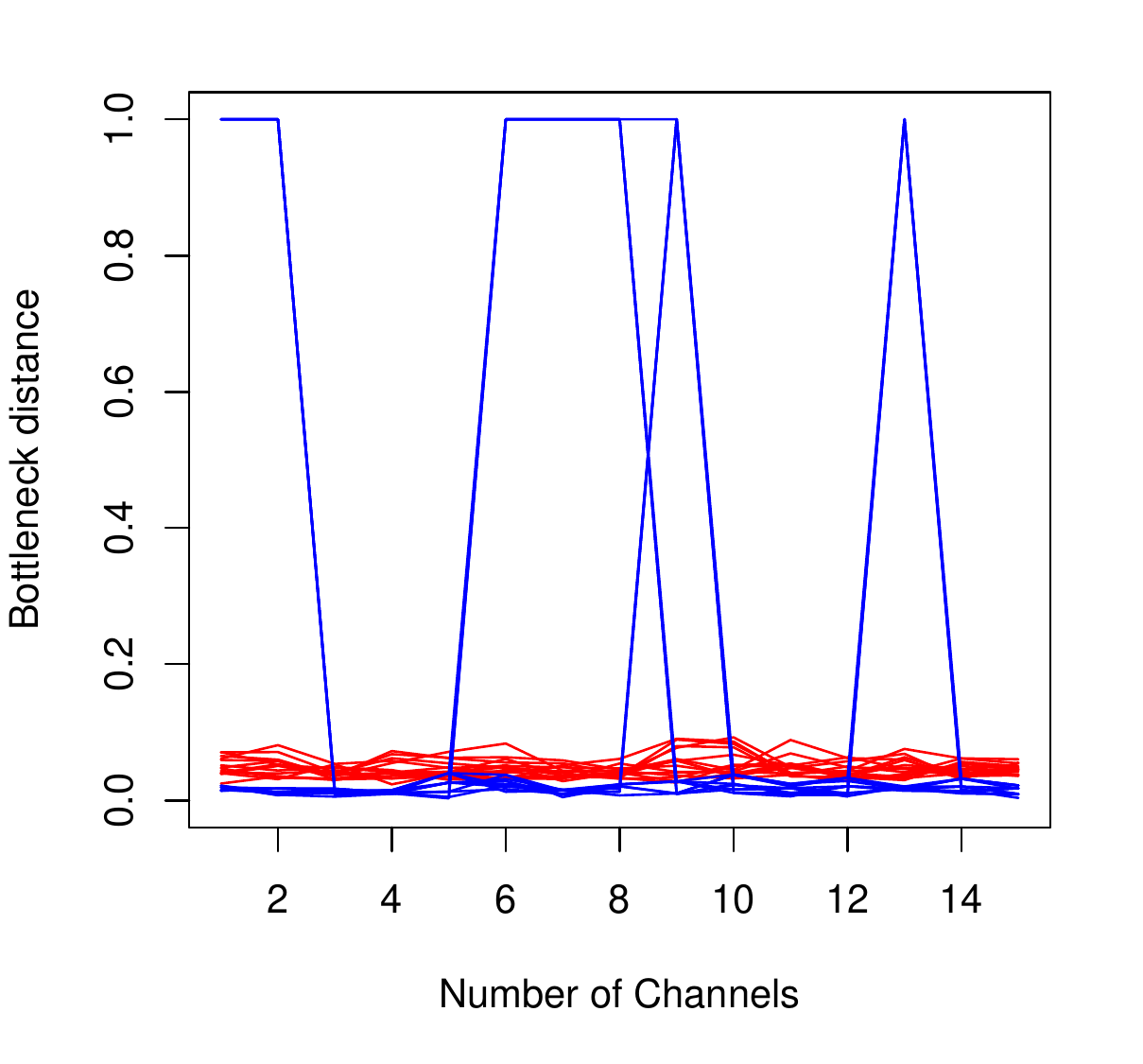}\\
(D) &  (E) & (C)\\
\includegraphics{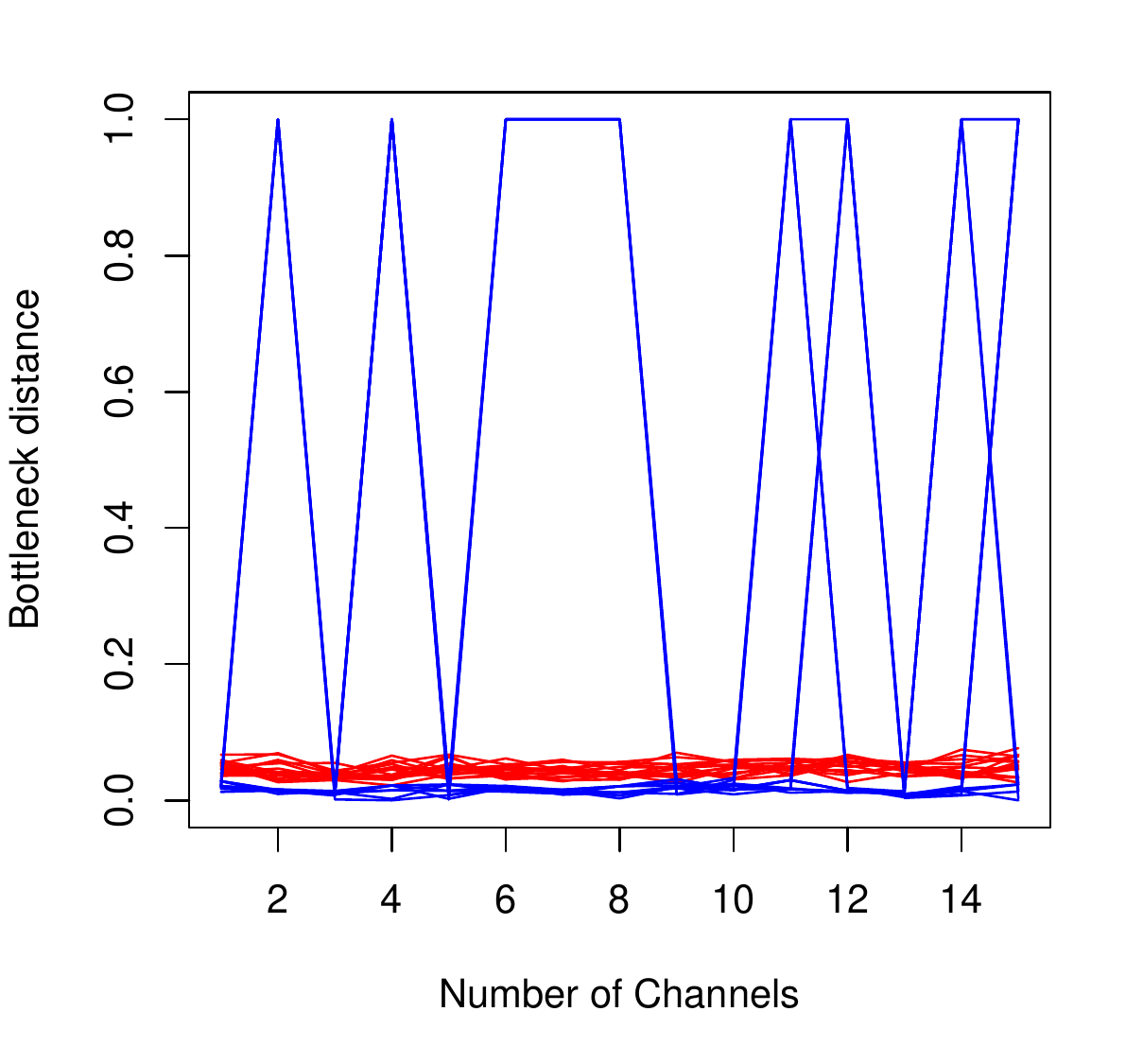} & \includegraphics{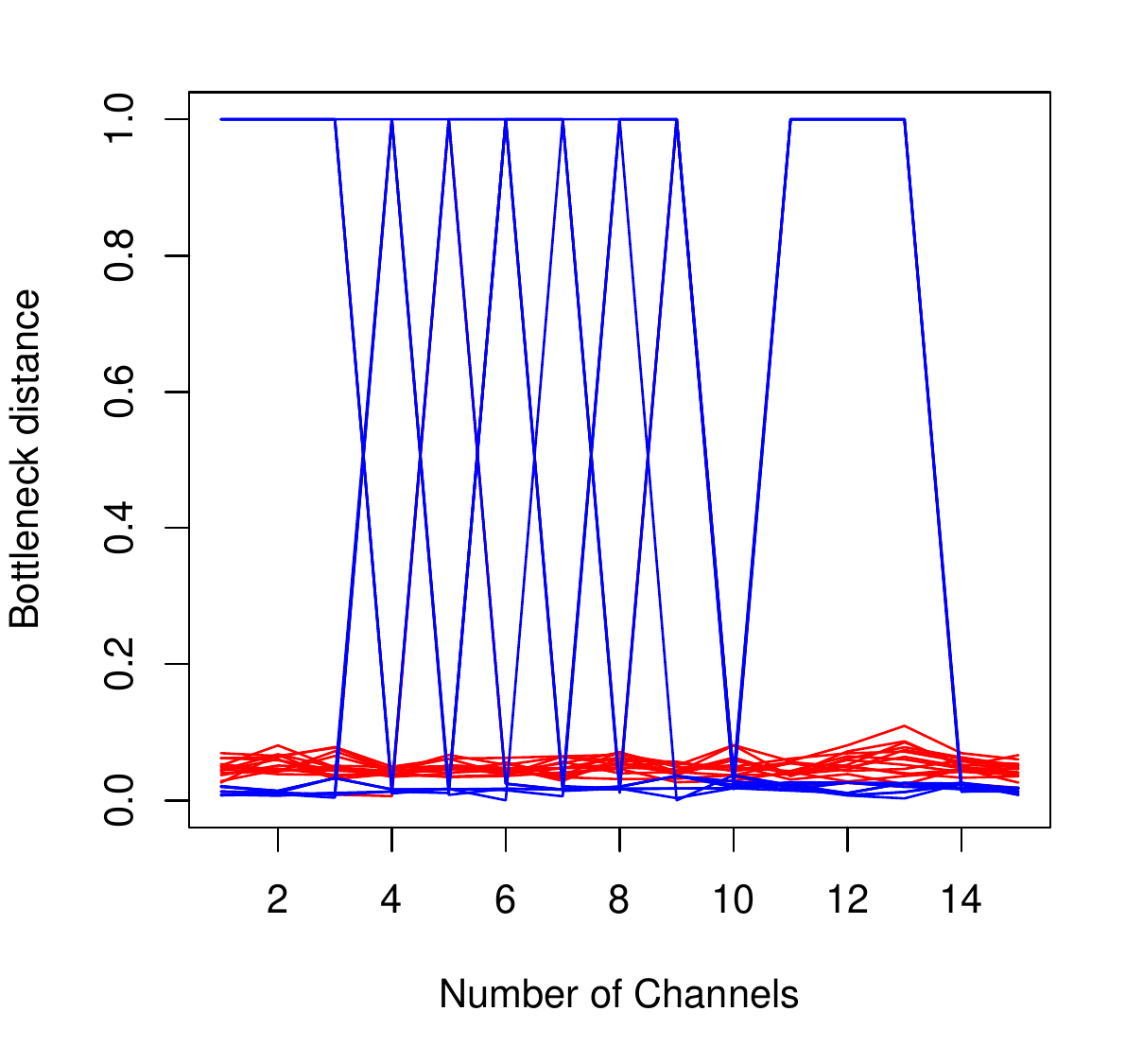} &

\end{tabular}
\end{minipage}}
\caption{Bottleneck distances between the persistent diagrams for 15 channels within each set A, B, C, D, and E on $\Hp_1$ and $\Hp_2$ for each of the methods introduced above.}
\label{fig7}
\end{figure}

\noindent (b) \noindent {Between set analysis} \\

\noindent To analyze the data of Bottleneck distances between sets, we need a summary statistics for each set from the data above. It is clear from Figure \ref{fig7} that the mean would not be a great summary statistics on $\Hp_1$, as there seems to be too many outliers. We will use the median instead and perform a pairwise Wilcoxon-Mann-Whitney test. Table \ref{table2} below shows the p.value on $\Hp_1$ and and $\Hp_2$. The take-away is that  the last row of the table suggests that set E is statistically topologically different from others on $\Hp_1$, at significance level $0.05$. In a way, this is  a confirmation of the results obtained in \cite{Kwessi2021} where set E (seizure) was already  shown to be statistically  different from other sets. 
\begin{table}[H]
\centering 

\begin{tabular}{|c|c|c|c|c|c|}
\hline
   \red{$\Hp_1$}/\blue{$\Hp_2$} & A & B & C & D & E\\
\hline
A &  &  \blue{0.1975936} & \blue{0.3049497} & \blue{0.2467548} &  \blue{0.7432987} \\
\hline
B & \red{0.3202554} &  &  \blue{0.3835209}& \blue{0.5066311} &  \blue{0.1707835} \\
\hline
C & \red{0.0832231} & \red{0.1322987} & &\blue{0.8356690}  & \blue{0.7088614} \\
\hline
D & \red{0.2012797} & \red{0.6292608} & \red{0.6292608} & & \blue{0.5067258} \\
\hline
E & \red{0.0049325} & \red{0.0157855} & \red{0.0157855} & \red{0.0114901} & \\
\hline
\end{tabular}
\caption{P.values of  Wilcoxon-Mann-Whitney tests between sets of median Bottleneck distances. }
\label{table2}
\end{table}
\section{Concluding remarks} \label{sec5}
In this paper, we have revisited the mathematical descriptions of six dimension reduction methods. We have given a brief introduction to  the very vast  topic  of persistent homology. We discussed how to apply persistent homology to data. In the presence of data (say in three dimension) obtained either by projecting the data from high dimension into smaller dimension (as in Takens) or by performing some sort of dimension reduction,  it is not always clear what we see or how different  one method is compared to another. From their mathematical description, they seem to represent different objects. Further, obtaining theoretically a clear discrimination procedure between these procedures seems a daunting if not an outright impossible task. Topology may offer a solution by looking at persistent artifacts through filtration. From Figure \ref{fig5}, it seems clear that the methods were different but Figure \ref{fig6} offers a different perspective. In the end, through calculation of Bottleneck distances and hypothesis tests, we can safely conclude that the methods are different topologically speaking, in that, the connected components, the tunnels and holes, the shells around cavities or  voids do not match perfectly. Since these methods are indiscriminately used in many applications, the message is that replication of results from one method to the next may not be guaranteed in the grand scheme of things. It does not however render them useless. In fact, our analysis is limited to one data set, meaning  that another data set may yield different conclusions. Further, due to cost in calculation, we were limited to only a handful of samples.  More, Wasserstein distance for $p<\infty$ are extremely costly  in time to calculate  on a regular computer. Even for $p=\infty$, the Bottleneck distance is also very costly in time to calculate, especially for $\Hp_0$. This explain why at some point, we did not provide the comparison on $\Hp_0$. We can infer from this analysis that topological persistent homologies do change dramatically at seizure, a finding already obtained in previous analyses, see \cite{Kwessi2021}. This suggests that looking at changes in homology landscapes could be a predictor of seizure. Given that some EEG epilepsy data are known to contain some deterministic chaos,  it be might worthwhile  to study  whether persistent homology can also be used  for better understanding of chaotic data in dynamical systems. 
\bibliography{EEGAnaly_Persistent_Homology}

\begin{thebibliography}{33}
\providecommand{\natexlab}[1]{#1}
\providecommand{\url}[1]{\texttt{#1}}
\expandafter\ifx\csname urlstyle\endcsname\relax
  \providecommand{\doi}[1]{doi: #1}\else
  \providecommand{\doi}{doi: \begingroup \urlstyle{rm}\Url}\fi

\bibitem[Belkin and Niyogi(2002)]{Belkin2002a}
M.~Belkin and P.~Niyogi.
\newblock Laplacian eigenmaps and spectral techniques for embedding and
  clustering.
\newblock In T.~G. Dietterich, S.~Becker, and Z.~Ghahramani, editors,
  \emph{Advances in Neural Information Processing Systems 14}, pages 585--591.
  Cambridge, MA: MIT Press, 2002.

\bibitem[Berry et~al.(2020)Berry, Chen, Cisewski-Kehe, and Fasy]{Berry2020}
E.~Berry, Y.-C. Chen, J.~Cisewski-Kehe, and B.T. Fasy.
\newblock Functional summaries of persistence diagrams.
\newblock \emph{Journal of Applied and Computational Topology}, 4:\penalty0
  211--262, 2020.

\bibitem[Bhattacharya et~al.(2015)Bhattacharya, Ghrist, and
  Kumar]{Bhattacharya2015}
S.~Bhattacharya, R.~Ghrist, and V.~Kumar.
\newblock Persistent homology for path planning in uncertain environments.
\newblock \emph{IEEETrans Robot}, 31:\penalty0 578--590, 2015.

\bibitem[Chan et~al.(2013)Chan, Carlsson, and Rabadan]{Chan2013}
J.~M. Chan, G.~Carlsson, and R.~Rabadan.
\newblock Topology of viral evolution.
\newblock \emph{PNAS}, 110\penalty0 (46):\penalty0 18566--18571, 2013.

\bibitem[Chung et~al.(2009)Chung, Bubenik, and Kim]{Chung2009}
M.~K. Chung, P.~Bubenik, and P.~T. Kim.
\newblock Persistence diagrams of cortical surface data.
\newblock In J.~L. Prince, D.~L. Pham, and K.~J. Myers, editors,
  \emph{Information processing in medical imaging. Lecture notes in computer
  science}, volume 5636 of \emph{Springer, Berlin}, pages 386--397, 2009.

\bibitem[Chung et~al.(2023)Chung, Ramos, Paiva, Prabharakaren, Nair, Meyerand,
  Hermann, Binder, and Struck]{Chung2023}
M.~K. Chung, C.~G. Ramos, J.~Paiva, F. B.~Mathis, V.~Prabharakaren, V.~A. Nair,
  E.~Meyerand, B.~P. Hermann, J.~R. Binder, and A.~F. Struck.
\newblock Unified topological inference for brainnetworks in temporal lobe
  epilepsy using thewasserstein distance.
\newblock \emph{ArXiv}, 2023.

\bibitem[de~Silva and Ghrist(2007)]{Silva2007}
V.~G. de~Silva and R.~Ghrist.
\newblock Coverage in sensor networks via persistent homology.
\newblock \emph{Algebraic Geom Topol}, 7:\penalty0 339--358, 2007.
\newblock \doi{DOI 10.1140/epjds/s13688-017-0109-5}.

\bibitem[Emmett et~al.(2016)Emmett, Schweinhart, and Rabad\'an]{Emmett2016}
K.~Emmett, N.~Schweinhart, and R.~Rabad\'an.
\newblock Multiscale topology of chromatin folding.
\newblock In \emph{Proceedings of the 9th EAIinternational conference on
  bio-inspired information and communications technologies}, BICT’15. ICST,
  pages 177--180, 2016.

\bibitem[Gameiro et~al.(2015)Gameiro, Hiraoka, Izumi, M., Mischaikow, and
  Nanda]{Gameiro2015}
M.~Gameiro, Y.~Hiraoka, S.~Izumi, Kram\'ar M., K.~Mischaikow, and K.~Nanda.
\newblock A topological measurement of proteincompressibility.
\newblock \emph{Jpn J Ind Appl Math}, 32:\penalty0 1--1--7, 2015.

\bibitem[Ghrist(2008)]{Ghrist2008}
R.~Ghrist.
\newblock Barcodes: The persistent topology of data.
\newblock \emph{Bull. Amer. Math. Soc.}, 45:\penalty0 61--75, 2008.

\bibitem[Giusti et~al.(2016)Giusti, Ghrist, and Bassett]{Giusti2016}
C.~Giusti, R.~Ghrist, and D.~Bassett.
\newblock Two’s company and three (or more) is a simplex.
\newblock \emph{J Comput Neurosci}, 41:\penalty0 1--14, 2016.

\bibitem[Guillemard et~al.(2013)Guillemard, Boche, Kutyniok, and
  Philipp]{Guillemard2013}
M.~Guillemard, H.~Boche, G.~Kutyniok, and F.~Philipp.
\newblock Persistence diagrams of cortical surface data.
\newblock In \emph{10th international conference on sampling theory and
  applications.}, pages 309--312, 2013.

\bibitem[Hinton and Roweis(2002)]{Hinton2002}
G.~E. Hinton and S.~Roweis.
\newblock Stochastic neighbor embedding.
\newblock In S.~Becker, S.~Thrun, and K.~Obermayer, editors, \emph{Advances in
  Neural Information Processing Systems}, volume~15. MIT Press, 2002.

\bibitem[Hotelling(1933)]{Hotelling}
H.~Hotelling.
\newblock Analysis of a complex of statistical variables into principal
  components.
\newblock \emph{Journal of Educationbal Psychology}, 24\penalty0
  (417--441):\penalty0 498--520, 1933.

\bibitem[Hyv\"arinen(1999)]{Hyvarinen1999}
A.~Hyv\"arinen.
\newblock Fast and robust fixed-point algorithms for independent component
  analysis.
\newblock \emph{IEEE Transactions on Neural Networks}, 13\penalty0
  (4--5):\penalty0 411--430, 1999.

\bibitem[Kwessi and Edwards(2021)]{Kwessi2021}
E.~Kwessi and L.~Edwards.
\newblock Analysis of eeg time series data using complex structurization.
\newblock \emph{Neural Computations}, 33\penalty0 (7), 2021.

\bibitem[Leibon et~al.(2008)Leibon, Pauls, Rockmore, and Savell]{Leibon2008}
G.~Leibon, S.~Pauls, D.~Rockmore, and R.~Savell.
\newblock Topological structures in the equities market network.
\newblock \emph{Proc Natl AcadSci USA}, 105:\penalty0 20589--20594, 2008.

\bibitem[Ma and Fu(2012)]{MaFu}
Y.~Ma and Y.~Fu.
\newblock \emph{Manifold Learning: Theory and Applications}.
\newblock CRC Press, 2012.

\bibitem[Maleti\'c et~al.(2016)Maleti\'c, Zhao, and Rajkovi\'c]{Maletic2016}
S.~Maleti\'c, Y.~Zhao, and M.~Rajkovi\'c.
\newblock Persistent topological features of dynamical systems.
\newblock \emph{Chaos}, 26\penalty0 (5):\penalty0 053105, 2016.
\newblock \doi{10.1063/1.4949472. PMID: 27249945}.

\bibitem[Mileyko et~al.(2011)Mileyko, Mukherjee, , and Harer]{Mileyko2011}
Y.~Mileyko, S.~Mukherjee, , and J.~Harer.
\newblock Probability measures on the space of persistence diagrams.
\newblock \emph{Inverse Problems}, 27:\penalty0 124007, 2011.

\bibitem[Naizait et~al.(2020)Naizait, Zhitnikov, and Lim]{Naizait2020}
G.~Naizait, A.~Zhitnikov, and L.-H. Lim.
\newblock Topology of deep neural networks.
\newblock \emph{Journal of Machine Learning Research}, 21:\penalty0 1--40,
  2020.

\bibitem[Otter et~al.(2017)Otter, Porter, Tillman, Grindrod, and
  Harrington]{Otter2017}
N.~Otter, M.~A. Porter, U.~Tillman, P.~Grindrod, and H.~A. Harrington.
\newblock A roadmap for the computationof persistent homology.
\newblock \emph{EPJ Data Science}, 6\penalty0 (7), 2017.
\newblock \doi{DOI 10.1140/epjds/s13688-017-0109-5}.

\bibitem[Pokorny et~al.(2016)Pokorny, Hawasly, and Ramamoorthy]{Pokorny2016}
F.T. Pokorny, M.~Hawasly, and S.~Ramamoorthy.
\newblock Topological trajectory classification with filtrations of
  simplicialcomplexes and persistent homology.
\newblock \emph{Int J Robot Res}, 35:\penalty0 204--223, 2016.

\bibitem[Rizvi et~al.(2017)Rizvi, Camara, Kandror, Roberts, Schieren, Maniatis,
  and Rabad\'an]{Rizvi2017}
A.~Rizvi, P.~Camara, E.~Kandror, T.~Roberts, I.~Schieren, T.~Maniatis, and
  R.~Rabad\'an.
\newblock Single-cell topological rna-seqanalysis reveals insights into
  cellular differentiation and development.
\newblock \emph{Nat Biotechnol}, 35:\penalty0 551--560, 2017.

\bibitem[Sizemore et~al.(2019)Sizemore, Phillips-Cremins, Ghrist, and
  Bassett]{Sizemore2019}
A.~E. Sizemore, J.~E. Phillips-Cremins, R.~Ghrist, and D.~S. Bassett.
\newblock The importance of the whole: Topological data analysis for the
  network neuroscientist.
\newblock \emph{Network Neuroscience}, 3:\penalty0 656–--673, 2019.

\bibitem[Takens(1981)]{Takens1981}
F.~Takens.
\newblock Detecting strange attractors in turbulence dynamical systems and
  turbulence (lecture notes in mathematics), vol. 898, 1981.

\bibitem[Taylor et~al.(2015)Taylor, Klimm, Harrington, Kram\'ar, Mischaikow,
  Porter, and Mucha]{Taylor2015}
D.~Taylor, F.~Klimm, H.~A. Harrington, M.~Kram\'ar, K.~Mischaikow, M.~A.
  Porter, and P.~J. Mucha.
\newblock Topological data analysis ofcontagion maps for examining spreading
  processes on networks.
\newblock \emph{Nat Commun}, 6, 2015.
\newblock Article ID 7723.

\bibitem[Tenenbaum et~al.(2000)Tenenbaum, de~Siva, and Langford]{Tenenbaum2000}
J.B. Tenenbaum, V.~de~Siva, and J.~C. Langford.
\newblock A global geometric frameworkfor nonlinear dimensionality reduction.
\newblock \emph{Science}, 290:\penalty0 2319--2323, 2000.

\bibitem[Torgerson(1952)]{Torgerson1952}
W.~S. Torgerson.
\newblock Multidimensional scaling: I. theory and method.
\newblock \emph{Psychometrika}, 17:\penalty0 410--419, 1952.

\bibitem[van~der Maaten and Hinton(2008)]{Maaten2008}
L.~J.~P. van~der Maaten and G.~E. Hinton.
\newblock Visualizing data using t-sne.
\newblock \emph{Journal of Machine Learning Research}, 9:\penalty0 2579--2605,
  2008.

\bibitem[Vasudevan et~al.(2013)Vasudevan, Ames, and Bajcsy]{Vasudevan2013}
R.~Vasudevan, A.~Ames, and R.~Bajcsy.
\newblock Persistent homology for automatic determination of human-data based
  costof bipedal walking.
\newblock \emph{Nonlinear Anal Hybrid Syst}, 7:\penalty0 101--115, 2013.

\bibitem[Whitney(1936)]{Whitney}
H.~Whitney.
\newblock Differentiable manifolds.
\newblock \emph{Ann. Math}, 37\penalty0 (3):\penalty0 645--680, 1936.

\bibitem[Xia and Wei(2014)]{Xia2014}
K.~Xia and G.-W. Wei.
\newblock Persistent homology analysis of protein structure, flexibility, and
  folding.
\newblock \emph{Int J Numer MethodsBiomed Eng}, 30:\penalty0 814--844, 2014.

\end{thebibliography}


 \end{document}